\ifdefined\XeTeXversion\else\pdfoutput=1\fi
\documentclass{article}
\usepackage[utf8]{inputenc}
\usepackage[T1]{fontenc}
\usepackage{lmodern}
\usepackage{arxiv}
\usepackage{authblk}
\usepackage[numbers,sort&compress]{natbib}

\usepackage{amssymb}
\usepackage{amsmath}
\usepackage{array}
\usepackage{booktabs}
\usepackage{graphicx}
\usepackage{url}
\usepackage{hyperref}
\hypersetup{hidelinks}
\usepackage{cleveref}
\usepackage{subcaption}
\crefname{section}{Section}{Sections}
\crefname{subsection}{Section}{Sections}
\crefname{figure}{Figure}{Figures}
\crefname{equation}{Eq.}{Eqs.}
\crefname{appendix}{Appendix}{Appendices}
\usepackage{float}
\usepackage{listings}
\usepackage{xcolor}
\lstdefinelanguage{Julia}{%
  morekeywords={function,end,return,for,in,if,else,elseif,while,using,const,let,do,true,false,nothing},
  morekeywords=[2]{ODEProblem,ShootingLayer,Tunable,Constant,CorleoneDynamicOptProblem,OptimizationProblem,
    solve,pack,default_vectorizer,OEDLayer,ContinuousMeasurement,ACriterion,Chain,Dense,setup,apply,collect},
  sensitive=true, morecomment=[l]{\#}, morestring=[b]"}
\definecolor{lstkeyword}{RGB}{0,90,160}
\definecolor{lstfunction}{RGB}{120,40,140}
\definecolor{lstcomment}{RGB}{90,110,90}
\definecolor{lststring}{RGB}{170,60,20}
\definecolor{lstbackground}{RGB}{247,247,247}
\crefname{lstlisting}{Listing}{Listings}
\crefname{table}{Table}{Tables}

\newcommand{\barr}{\begin{array}}
\newcommand{\earr}{\end{array}}
\newcommand{\eqleq}{\mathrel{\substack{=\\[-0.2ex] \leq}}}

\newcommand{\papertitle}{Corleone.jl: A modular shooting framework for dynamic optimization and scientific machine learning in Julia}
\newcommand{\paperauthors}{Carl Julius Martensen, Christoph Plate, Christian Kirches, Andreas Potschka, Christopher Vincent Rackauckas, Sebastian Sager}
\newcommand{\paperkeywords}{dynamic optimization \kwsep optimal control \kwsep direct multiple shooting \kwsep scientific machine learning \kwsep universal differential equations \kwsep Julia}

\hypersetup{
  pdftitle={\papertitle},
  pdfauthor={\paperauthors},
  pdfkeywords={dynamic optimization, optimal control, direct multiple shooting, scientific machine learning, universal differential equations, Julia}
}

\title{\papertitle}

\author[1]{Carl~Julius~Martensen}
\author[1,2]{Christoph~Plate}
\author[3]{Christian~Kirches}
\author[4]{Andreas~Potschka}
\author[5]{Christopher~Vincent~Rackauckas}
\author[1,2]{Sebastian~Sager}
\affil[1]{Otto von Guericke University Magdeburg\protect\\ Universit\"atsplatz~2, 39106 Magdeburg, Germany}
\affil[2]{\protect\rule{0pt}{3ex}Max Planck Institute for Dynamics of Complex Technical Systems\protect\\ Sandtorstra\ss e~1, 39106 Magdeburg, Germany}
\affil[3]{\protect\rule{0pt}{3ex}Braunschweig University of Technology\protect\\ Universit\"atsplatz~2, 38106 Braunschweig, Germany}
\affil[4]{\protect\rule{0pt}{3ex}Clausthal University of Technology\protect\\ Erzstra\ss e~1, 38678 Clausthal-Zellerfeld, Germany}
\affil[5]{\protect\rule{0pt}{3ex}Massachusetts Institute of Technology\protect\\ 77 Massachusetts Avenue, Cambridge, MA 02139, USA}
\affil[ ]{\protect\rule{0pt}{3ex}Correspondence: \texttt{carl.martensen@ovgu.de}}

\newcommand{\kwsep}{\and}
\newcommand{\contributionsheading}{Author contributions}
\newcommand{\competingheading}{Competing interests}
\newcommand{\aiheading}{Use of generative AI}
\newcommand{\thisarticle}{this article}

\begin{document}
\renewcommand{\labelenumii}{\arabic{enumi}.\arabic{enumii}}
\maketitle

\begin{abstract}
Corleone.jl is a Julia package for dynamic optimization problems with mechanistic, data-driven, or hybrid models. It represents the piecewise simulation of a dynamical system as a Lux.jl layer, so that single and multiple shooting, trainable model components, and differential equation solvers of the scientific machine learning (SciML) ecosystem are combined within one differentiable framework. On a technical level, Corleone.jl acts as a transcription layer that provides objectives, constraints, derivatives, and sparsity patterns to (mixed-integer) nonlinear programming solvers. We illustrate several use cases---optimal control, including mixed-integer nonlinear optimal control, the training of universal differential equations, and optimal experimental design---and outline available extensions to inverse problems and dynamic games.

\end{abstract}


\section{Motivation and significance} \label{sec:motivation}

Dynamical systems describe evolving processes across chemical~\cite{srinivasanDynamicOptimizationBatch2003} and mechanical engineering~\cite{hrovatDevelopmentModelPredictive2012}, pharmacometrics~\cite{mooreHowMathematicallyOptimize2018}, biology~\cite{raueLessonsLearnedQuantitative2013}, robotics~\cite{wensingOptimizationbasedControlDynamic2024}, and energy grids~\cite{huModelPredictiveControl2021}. Using these models for scientific inference or industrial applications often requires solving an optimization problem constrained by the dynamics. Examples include estimating unknown parameters or initial conditions, designing optimal experiments, reconstructing unobserved states, and computing optimal time-dependent inputs. Such dynamic optimization problems are commonly solved by direct methods, which discretize the controls and transcribe the problem into a nonlinear program (NLP)~\cite{bettsPracticalMethodsOptimal2010, bieglerNonlinearProgrammingConcepts2010}; see \cite{sassObscuredTerminologyHidden2025} for a tutorial on the terminology of these methods. Direct multiple shooting \cite{bockplitt_1984_jul_multipleshootingalgorithm} is one such transcription, and structure-exploiting NLP algorithms \cite{leineweberEfficientMultipleShooting2003, raoApplicationInteriorpointMethods1998, frisonHPIPMHighperformanceQuadratic2020, jankaetal_2016_dec_sr1bfgssqp, wittmannsager___blocksqp2exploiting} make use of the block structure of the resulting NLPs. Depending on the task at hand, specialized algorithmic approaches \cite{burger_pycombina_2020} can greatly reduce the solution time. Tailored formulations, e.g., for optimal experimental design \cite{bauerNumericalMethodsOptimum2000, franceschiniModelbasedDesignExperiments2008, sager_2013_jan_samplingdecisionsoptimum}, combine practical requirements with sound theoretical foundations and have led to specialized software \cite{wangPyomoDOEOpensourcePackage2022, martensenDynamicOEDjlJuliaPackage2024}.
Scientific machine learning (SciML) creates new opportunities \cite{karniadakisPhysicsinformedMachineLearning2021, rackauckasetal_2021_nov_universaldifferentialequations} and challenges \cite{plateetal_2026_mar_optimalexperimentaldesign} for classical workflows and practitioners alike, and requires software that combines trainable components with mechanistic models and numerical routines.

Existing tools take different approaches to address these challenges. Symbolic frameworks such as CasADi \cite{anderssonCasADiSoftwareFramework2019} let users write collocation or shooting transcriptions themselves, whereas modeling environments such as Pyomo, InfiniteOpt.jl, OptimalControl.jl, and ModelingToolkit.jl default to direct transcription by collocation or fixed-step Runge--Kutta schemes \cite{nicholsonPyomodaeModelingAutomatic2018, pulsipherUnifyingModelingAbstraction2022, caillauOptimalControljlJuliaPackage2026, maModelingToolkitComposableGraph2021}. In all cases, the modeling interface, the differentiation, and the integration backend determine which existing models can be used, and machine-learning components typically enter through a separate bridge \cite{salzmannLearningCasADiDatadriven2024, lahrL4acadosLearningbasedModels2024}. Practitioners whose models already use adaptive integrators, sensitivity algorithms, or neural-network layers in a differentiable ecosystem therefore have to re-express or wrap them before they can start to optimize.

\paragraph{Requirements} Industrial applications, in particular in the chemical and automotive sectors, and scientific ones, in particular in the life sciences, motivate the following list of requirements for dynamic optimization software:
\begin{description}\setlength{\itemsep}{0pt}\setlength{\parsep}{0pt}\setlength{\topsep}{0pt}
\item[(R1)] Active learning and optimal experimental design.
\item[(R2)] Calibration of models, i.e., parameter estimation and/or training.
\item[(R3)] Complexity reduction, including model reduction, parameter-space reduction, structural preprocessing, and condensing.
\item[(R4)] Mixed-integer optimization.
\item[(R5)] Multistage problems, i.e., discontinuities or even changes of the model in time within a single optimization problem.
\item[(R6)] Use of hybrid models, i.e., mechanistic models with trainable components such as neural networks.
\item[(R7)] Reuse of existing models and industrial standards.
\item[(R8)] Compiled code for runtime performance.
\item[(R9)] Parallelization.
\item[(R10)] Sustainability and modularity within an active software ecosystem, in particular through interfaces for structure exploitation.
\item[(R11)] Availability under a permissive license approved by the Open Source Initiative (OSI).
\end{description}
 R1--R6 concern the problem classes that the software has to handle, and R7--R11 its integration into research and industrial workflows.
\Cref{tab:requirements} relates these requirements to established packages for optimal control, model predictive control, algebraic modeling, model calibration, and optimal experimental design, several of which cover a broad subset of them.

\begin{table}[h!!!]
\centering
\small
\setlength{\tabcolsep}{3pt}
\begin{tabular}{lccccccccccc}
\hline
Software & R1 & R2 & R3 & R4 & R5 & R6 & R7 & R8 & R9 & R10 & R11 \\
\hline
acados~\cite{verschuerenAcadosModularOpensource2022, lahrL4acadosLearningbasedModels2024} &  & \checkmark & \checkmark &  & \checkmark & (\checkmark) &  & \checkmark & \checkmark & \checkmark & \checkmark \\
AMICI/pyPESTO~\cite{frohlichAMICIHighperformanceSensitivity2021, schaltePyPESTOModularScalable2023, schmiesterPEtabInteroperableSpecification2021} &  & \checkmark & (\checkmark) &  & \checkmark & \checkmark & \checkmark & \checkmark & \checkmark & \checkmark & \checkmark \\
AMIGO2~\cite{balsacantoAMIGO2ToolboxDynamic2016} & \checkmark & \checkmark &  &  &  &  & (\checkmark) & \checkmark &  &  &  \\
BOCOP~\cite{bonnansBocopCollectionExamples2012} &  & \checkmark &  &  &  &  &  & \checkmark &  &  & *  \\
CasADi~\cite{anderssonCasADiSoftwareFramework2019, salzmannLearningCasADiDatadriven2024} &  & \checkmark & (\checkmark) & \checkmark & \checkmark & (\checkmark) & \checkmark & \checkmark & \checkmark & \checkmark &  * \\
do-mpc~\cite{fiedlerDompcFAIRNonlinear2023} &  & \checkmark &  &  &  & \checkmark &  & (\checkmark) &  &  &  * \\
Dymos~\cite{falckDymosPythonPackage2021} &  & (\checkmark) &  &  & \checkmark &  & (\checkmark) &  & \checkmark & \checkmark & \checkmark \\
DyOS~\cite{caspariDyOSFrameworkOptimization2019} &  &  &  & \checkmark & \checkmark &  & \checkmark & \checkmark &  &  &  \\
GEKKO~\cite{bealGEKKOOptimizationSuite2018} &  & \checkmark & \checkmark & \checkmark & \checkmark & \checkmark &  & \checkmark &  &  & \checkmark \\
GPOPS-II~\cite{pattersonGPOPSIIMATLABSoftware2014, agamawiCGPOPSOptimalControl2020} &  &  &  &  & \checkmark &  &  & (\checkmark) &  &  &  \\
HILO-MPC~\cite{pohlodekFlexibleDevelopmentEvaluation2025} &  & \checkmark & (\checkmark) &  &  & \checkmark &  & (\checkmark) &  &  &  * \\
ICLOCS2~\cite{nieICLOCS2TryThis2018} &  & (\checkmark) &  &  & \checkmark &  &  &  &  &  & \checkmark \\
InfiniteOpt.jl~\cite{pulsipherUnifyingModelingAbstraction2022} &  & \checkmark &  & \checkmark &  & (\checkmark) & (\checkmark) & \checkmark & \checkmark & \checkmark & \checkmark \\
MUSCOD-II~\cite{hoffmannMUSCODIIUserManual2024} &  & \checkmark & \checkmark & \checkmark & \checkmark &  &  & \checkmark &  &  &  \\
NeuroMANCER~\cite{drgonaDifferentiablePredictiveControl2022} &  & \checkmark &  & \checkmark &  & \checkmark &  &  & \checkmark & \checkmark & \checkmark \\
OptimalControl.jl~\cite{caillauOptimalControljlJuliaPackage2026} &  & (\checkmark) &  &  &  &  &  & \checkmark & \checkmark & \checkmark & \checkmark \\
PEtab.jl~\cite{perssonPEtabjlAdvancingEfficiency2025} &  & \checkmark &  &  & \checkmark & \checkmark & \checkmark & \checkmark &  & \checkmark & \checkmark \\
PSOPT~\cite{becerraSolvingComplexOptimal2010} &  & \checkmark &  & \checkmark & \checkmark &  &  & \checkmark &  &  &  * \\
PyOED~\cite{chowdharyPyOEDExtensibleSuite2024} & \checkmark & \checkmark & \checkmark & (\checkmark) &  &  &  &  &  &  & \checkmark \\
Pyomo~\cite{hartPyomoModelingSolving2011, nicholsonPyomodaeModelingAutomatic2018, wangPyomoDOEOpensourcePackage2022} & \checkmark & \checkmark & \checkmark & \checkmark &  & (\checkmark) &  &  & \checkmark & \checkmark & \checkmark \\
VPLAN~\cite{koerkelNumerischeMethodenOptimale2002} & \checkmark & \checkmark &  &  &  &  &  & \checkmark &  &  &  \\
\hline
Corleone.jl (this paper) & \checkmark & \checkmark & (\checkmark) & (\checkmark) & \checkmark & \checkmark & (\checkmark) & \checkmark & \checkmark & \checkmark & \checkmark \\
\hline
\end{tabular}
\caption{Requirements, assessed from primary references and documentation as of November 2026. A check mark indicates documented support in the package or a maintained ecosystem extension. Parentheses indicate restricted support: calibration via generic parameter optimization (R2), solver-specific, linear-only, or task-specific reduction (R3), binary sensor selection or integrality passed to generic MINLP solvers without dedicated mixed-integer optimal control methods (R4), an external learning bridge (R6), model reuse without demonstrated standards support (R7), or restricted code generation (R8). R8 requires model and derivative code compiled to machine code, just in time or by code generation; R10 requires public interfaces that pass problem structure, such as sparsity patterns, to solvers. Blank entries mean that support was not established in the consulted sources. OSI-approved weak-copyleft licenses are marked with an asterisk (R11); they permit the use within proprietary software but require that modifications of the package itself are published under the same license.}
\label{tab:requirements}
\end{table}

\paragraph{Contribution} Corleone.jl (COmplexity Reduction and Lifting Enhanced Optimization of Nonlinear processes and Experiments) addresses these requirements jointly within one framework and one language, Julia~\cite{bezansonJuliaFreshApproach2017}. It represents dynamic simulations as layers of the deep-learning library Lux.jl~\cite{palLuxExplicitParameterization2023}, in which a layer is a function of its input, its trainable parameters, and a non-trainable state. It connects the differential equation solvers~\cite{rackauckasDifferentialEquationsjlPerformantFeaturerich2017} of the SciML ecosystem, a collection of Julia packages for differential equations, sensitivity analysis, and optimization, to trajectory-based objectives, constraints, and numerical optimization~\cite{dixitOptimizationjlUnifiedOptimization2023}. On the technical level, it is a transcription layer between a dynamic optimization problem and the finite-dimensional problem presented to an NLP solver; on the user level, it provides a common foundation for optimizing controls, model parameters, and differential states at the shooting nodes, and for constructing specialized extensions.
With respect to \cref{tab:requirements}, optimal experimental design (OED, R1) is provided by CorleoneOED.jl. Complexity reduction (R3) is currently restricted to OED: CorleoneOED.jl reduces the parameter space of the information matrix by a singular value decomposition, which makes OED for universal differential equations tractable~\cite{plateetal_2026_mar_optimalexperimentaldesign}. Hybrid models (R6) are native: trainable components of the dynamics enter as tunable parameters of the wrapped differential equation problem, and the simulation itself is a Lux.jl layer that can be composed with other layers. Existing models (R7) can be reused as SciML problems or as systems of ModelingToolkit.jl~\cite{maModelingToolkitComposableGraph2021}; industrial standards such as the Functional Mock-up Interface (FMI)~\cite{blochwitzFunctionalMockupInterface2011} are to be supported through ModelingToolkit.jl. For structure-exploiting solvers, Corleone.jl provides the sparsity patterns of constraint Jacobians and Lagrangian Hessians (R10). The contribution is the integration of these capabilities with Lux.jl layers and SciML simulation and differentiation workflows, so that a model written once for simulation or training can be optimized with the same integrator, sensitivity algorithm, and trainable components.

\section{Software description} \label{sec:software}

\subsection{Problem class} \label{sec:problemclass}

Corleone.jl targets dynamic optimization problems with mechanistic, data-driven, or hybrid models. In hybrid models, known physical or biological mechanisms are combined with trainable components such as neural networks~\cite{psichogiosHybridNeuralNetworkfirst1992, schweidtmannReviewPerspectiveHybrid2024}. Neural ordinary differential equations replace the entire right-hand side by a neural network~\cite{chenNeuralOrdinaryDifferential2018, kidgerNeuralDifferentialEquations2021}, universal differential equations (UDEs) replace only parts of an otherwise mechanistic model~\cite{rackauckasetal_2021_nov_universaldifferentialequations}. From the perspective of optimization, the weights of the trainable components are additional time-invariant parameters.

We consider multistage problems of the form
\begin{equation}
\barr{rcll}
\multicolumn{3}{l}{\displaystyle \min_{u_i(\cdot),\, q_i,\, c,\, x_i(\cdot)} \quad \sum_{i=1}^{N} \phi_i[x_i,u_i,q_i,c]} & \\[2ex]
\text{s.t. } \quad 0 & = & F_i\bigl(\dot{x}_i(t), x_i(t), u_i(t), q_i, c\bigr), & t \in [s_{i-1}, s_i],\\[0.5ex]
x_{i}(s_{i-1}) & = & F_{\mathrm{tr},i}\bigl(x_{i-1}(s_{i-1}), q_i, c\bigr), & \\[1.3ex] 
0 & \leq & h_i\bigl(x_i(t), u_i(t), q_i, c\bigr), & t \in [s_{i-1}, s_i],\\[1.1ex] 
0 & \eqleq & \displaystyle  r\bigl( x_{i(1)}(\tau_{1}), q_{i(1)}, \dots, x_{i(n_\tau)}(\tau_{n_\tau}), q_{i(n_\tau)}, c\bigr), & \\[2ex]
u_i(t) & \in & \mathcal{U}_i, \quad q_i \in \mathcal{Q}_i, \quad c \in \mathcal{C}, & t \in [s_{i-1}, s_i],
\earr
\label{eq:problemclass}
\end{equation}
where all constraints hold for the stages $i = 1, \dots, N$, on a fixed stage grid $s_0 < s_1 < \dots < s_N$. On stage $i$, the state vector is $x_i$ (note that dimensions may change from stage $i-1$ to $i$), $u_i$ the controls, and $q_i$ stage-wise time-invariant parameters; $c$ is the vector of global trainable parameters and weights. The admissible sets $\mathcal{U}_i = \{u: [s_{i-1}, s_i] \to \mathbb{R}^{n_{u_i}} \mid u_k: [s_{i-1}, s_i] \to \mathbb{Z},\ k \in I_{u_i}\}$, $\mathcal{Q}_i = \{q \in \mathbb{R}^{n_{q_i}} \mid q_k \in \mathbb{Z},\ k \in I_{q_i}\}$, and $\mathcal{C} = \{c \in \mathbb{R}^{n_c} \mid c_k \in \mathbb{Z},\ k \in I_c\}$ restrict the components with indices in $I_{u_i}$, $I_{q_i}$, and $I_c$ to integer values. These sets may contain additional combinatorial constraints. Corleone.jl thus handles mixed-integer nonlinear optimal control problems: it transcribes them into mixed-integer nonlinear programs (MINLP), which are passed to specialized solvers through Optimization.jl, or solves their relaxations for use in top-level algorithms (\cref{sec:example-control}). The stage objectives $\phi_i$ may comprise Mayer terms, which evaluate the states at the end of the stage, and Lagrange terms, which integrate a function of the states, parameters, and controls over the stage. Least-squares objectives of parameter estimation fit this form as well: continuous least-squares terms are Lagrange terms, and discrete ones evaluate the states at measurement times within the stage. The implicit model equations $F_i$ comprise ordinary differential equations and those differential-algebraic equations that are supported by SciMLSensitivity.jl and automatic differentiation, the transition functions $F_{\mathrm{tr},i}$ map the final state of stage $i-1$ to the initial state of stage $i$. A special case are initial values $x_0$, which are either given or tunable (and then formally part of $c$ in parameter estimation or of $q_1$ in control). We call a quantity tunable if it is an optimization variable of problem~\eqref{eq:problemclass} rather than a given constant; the term follows the annotations \texttt{Tunable} and \texttt{Constant} of Corleone.jl (\cref{sec:functionalities}). The mixed state-control constraints $h_i$ restrict states and controls, and the point constraints $r$ act on the states at the measurement or constraint times $\tau_j \in \{\tau_1, \dots, \tau_{n_\tau}\}$. They can link states at different times, e.g., for periodicity constraints. 
The times $\tau_j$ form a grid of their own, which may, but need not, contain the stage grid points, and $i(j)$ denotes the stage that contains $\tau_j$. After discretization, the mixed state-control constraints $h_i$ are relaxed to this grid as well. The symbol $\eqleq$ indicates either equality or inequality constraints. Parameter estimation, the training of hybrid models, optimal control, and optimal experimental design are special cases of problem~\eqref{eq:problemclass}.

Corleone.jl solves problem~\eqref{eq:problemclass} by a first-discretize-then-optimize approach. Each stage is simulated by solving initial value problems on a sequence of intervals. At selected interval boundaries, the state is reset to a new value instead of being carried over from the preceding interval. If such a reset value is tunable, the difference between the reset value and the state at the end of the preceding interval becomes a residual that the optimization drives to zero. Single and multiple shooting~\cite{bockplitt_1984_jul_multipleshootingalgorithm} are therefore special cases of the same construction: single shooting has no interior resets, whereas multiple shooting places tunable resets for all states at the shooting nodes, so that the residuals become the continuity constraints of the optimization problem. Multiple shooting has also been used to train neural ODEs~\cite{massaroliDifferentiableMultipleShooting2021, turanMultipleShootingTraining2022}.

\begin{figure}[t]
    \centering
    \includegraphics[width=\linewidth]{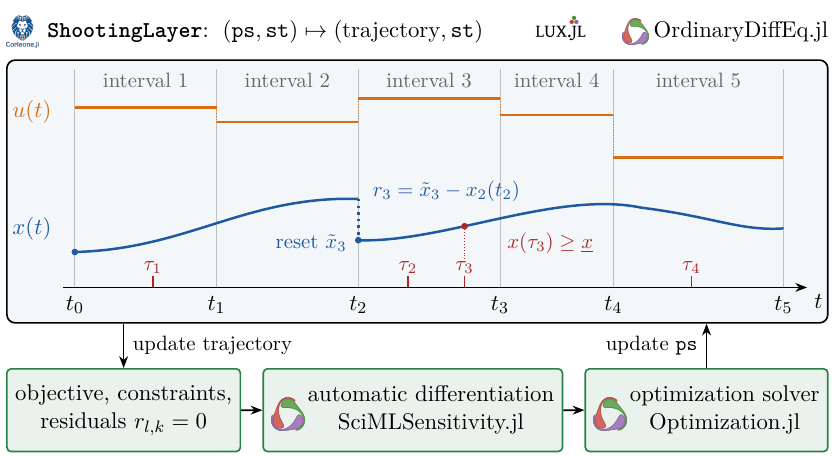}
    \caption{Forward simulation of one stage in Corleone.jl, showing a \texttt{ShootingLayer} and the surrounding optimization loop. The interval grid $t_0, \dots, t_5$ refines the stage $[s_{i-1}, s_i]$ of problem~\eqref{eq:problemclass}, i.e., $t_0 = s_{i-1}$ and $t_5 = s_i$. The layer is built on Lux.jl and integrates each interval with OrdinaryDiffEq.jl, as indicated by the logos at the top right. The dynamics of a single state $x$ depend on a piecewise constant control $u$. At $t_2$, a tunable value $\tilde{x}_3$ resets the trajectory, which yields the shooting residual $r_3$. Constraints are evaluated on a grid of times $\tau_j$ of its own, which may, but need not, contain the interval grid points; at $\tau_3$, a constraint bounds the state from below.
    The trajectory enters the objective and constraints of the optimization problem, whose derivatives are computed by automatic differentiation. AI generated.}
    \label{fig:sequence}
\end{figure}

Within a stage, each control is discretized as piecewise constant on its own control grid, and each state component $x_{\cdot,k}$ may be reset at its own reset times $T_{X_k}$. The interval grid $s_{i-1} = t_0 < t_1 < \dots < t_M = s_i$ of stage $i$ contains all control grid points and reset times of the stage, so that every control is constant on each interval, $u(t) \equiv u_l$ for $t \in [t_{l-1}, t_l]$; here and in the following, we drop the stage index $i$. On interval $l$, Corleone.jl solves the initial value problem
\begin{equation}
  \begin{aligned}
    0 &= F\bigl(\dot{x}_l(t), x_l(t), u_l, q, c\bigr), \qquad t \in [t_{l-1}, t_l],\\
    x_{l,k}(t_{l-1}) &=
    \begin{cases}
      \tilde{x}_{l,k} & \text{if } t_{l-1} \in T_{X_k},\\
      x_{l-1,k}(t_{l-1}) & \text{otherwise,}
    \end{cases}
  \end{aligned}
  \label{eq:ivp}
\end{equation}
for each state component $k$, where the reset value $\tilde{x}_{l,k}$ may be tunable or fixed. The first grid point $t_0$ belongs to every $T_{X_k}$; on the first stage, the values $\tilde{x}_{1,k}$ are the initial values $x_0$, and on later stages, they are given by the transition function $F_{\mathrm{tr},i}$ (\cref{sec:functionalities}). Without a reset, the state is carried over to the next interval. For differential-algebraic equations, the algebraic components of the state are made consistent with the reset values by the initialization of the integrator. Whenever a reset value $\tilde{x}_{l,k}$ is tunable, the trajectory stores the residual $r_{l,k} = \tilde{x}_{l,k} - x_{l-1,k}(t_{l-1})$. \Cref{fig:sequence} illustrates this construction for one stage.

\subsection{Software architecture} \label{sec:architecture}

Corleone.jl is built around the sequential solution of the initial value problems~\eqref{eq:ivp} on all intervals of a stage. Each interval defines the active control values and a possible reset of the initial conditions. This procedure is implemented as a Lux.jl layer~\cite{palLuxExplicitParameterization2023}, the \texttt{ShootingLayer}, which simulates one stage of problem~\eqref{eq:problemclass} and consists of sublayers for the initial conditions and the parameters, a numerical integrator, and an ensemble algorithm, which solves independent blocks of intervals serially, in threads, or in distributed processes. Such blocks arise from full-state resets, i.e., resets of all states at the same time; within a block, and when reverse-mode rules record the solves, the intervals are solved sequentially. Each interval is solved with a compatible integrator from DifferentialEquations.jl~\cite{rackauckasDifferentialEquationsjlPerformantFeaturerich2017}, including adaptive and stiff methods. Following the Lux.jl convention, the layer maps all tunable values, denoted by \texttt{ps} (parameters), such as time-invariant parameters, control values, and initial conditions, together with the non-differentiable layer state, denoted by \texttt{st} (state), such as the time grid, to a trajectory and the updated layer state, $(\texttt{ps}, \texttt{st}) \mapsto (\text{trajectory}, \texttt{st})$. A trajectory contains the solutions on all intervals, the residuals of all tunable resets, and the active control values at every time point. Because the simulation is a layer, specialized problems are built by wrapping or composing it rather than by reimplementing it: the \texttt{OEDLayer} of CorleoneOED.jl wraps a \texttt{ShootingLayer} of the model augmented with its sensitivities, and a multistage problem is a \texttt{Stages} container of shooting layers linked by \texttt{StageTransition} layers. Each of them exposes the same parameter, state, and trajectory interface to the optimization bridge. Conversely, a \texttt{ShootingLayer} can be part of a larger Lux.jl model, e.g., as the latent dynamics between the encoder and the decoder of a dynamic autoencoder.

\begin{figure}[t]
    \centering
    \includegraphics[width=\linewidth]{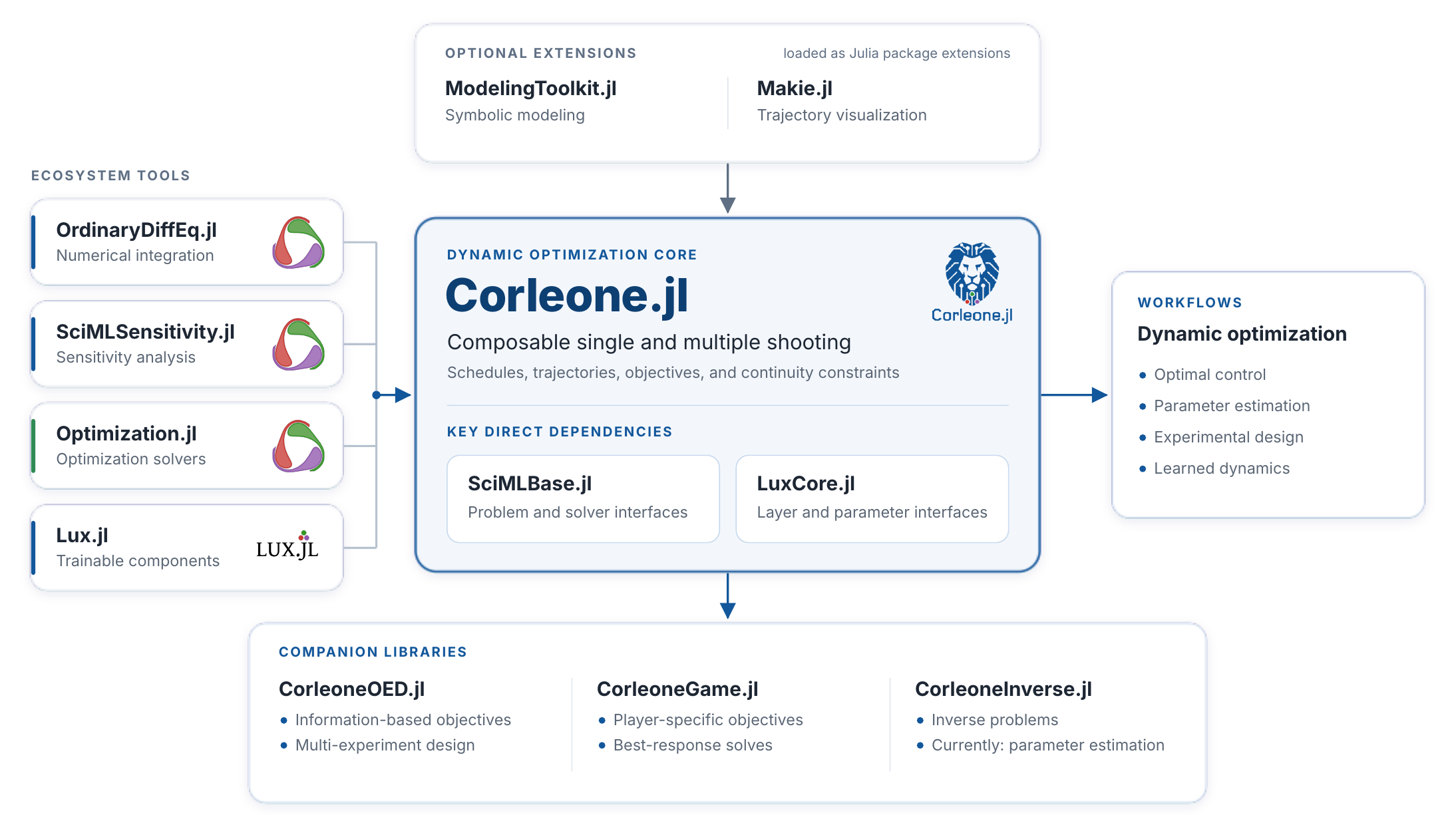}
    \caption{Architecture of Corleone.jl. Arrows indicate integrated functionality. The SciML ecosystem packages (left) supply integration, sensitivity analysis, optimization, and trainable components; the core (center) composes them into single and multiple shooting via the shared interfaces of SciMLBase.jl. Optional extensions (top) and companion libraries (bottom) implement specialized functionality on top of the core. AI generated.}
    \label{fig:architecture}
\end{figure}

Because the entire sequence is written in Julia, gradients with respect to all tunable values can be computed by automatic differentiation through the whole simulation; the test suite checks forward-mode gradients with ForwardDiff.jl and reverse-mode gradients with ReverseDiff.jl and with Zygote.jl and SciMLSensitivity.jl~\cite{revelsForwardModeAutomaticDifferentiation2016, innesDonUnrollAdjoint2018, sapienza_2026_diffeqad}. \Cref{fig:architecture} places the \texttt{ShootingLayer} within its dependencies, extensions, and companion libraries.

\subsection{Software functionalities} \label{sec:functionalities}

In addition to the \texttt{ShootingLayer}, Corleone.jl provides several features that simplify its use. Annotations (\texttt{Tunable} and \texttt{Constant}) declare for each state and parameter whether it is tunable, together with its default values or initializers, bounds, and integrality. Each annotation is paired with a time grid, the control grid of a control or the reset times $T_{X_k}$ of a state in \cref{eq:ivp}, and passed to the \texttt{ShootingLayer}, which uses both to initialize the parameters and initial conditions. Through SymbolicIndexingInterface.jl, annotations can refer to states and parameters by name rather than by position.

Objectives and constraints can be supplied either as plain Julia functions of the trajectory or as Julia expressions that refer to states and parameters by name and time. These are quoted Julia expressions \texttt{:(...)} rather than a separate modeling language. For example, the expression \texttt{:(x(12.0) - x(0.0))} as an equality constraint requires the state $x$ at the final time $t = 12$ to return to its initial value; this periodicity condition is a coupled multi-point constraint $r^c_j$ in problem~\eqref{eq:problemclass}. Controls and parameters are addressed in the same way. Corleone.jl collects all requested time points, adds them to the evaluation grid, and replaces each timed call by an index into the saved trajectory. The resulting expressions are compiled with \texttt{Core.eval}, or with RuntimeGeneratedFunctions.jl when it is loaded.

Code for specialized applications is loaded conditionally through Julia's package extension system. Examples are symbolic modeling with ModelingToolkit.jl, reverse-mode differentiation rules for Zygote.jl, and sparse Jacobians and Hessians that exploit the time structure of the problem with DifferentiationInterface.jl~\cite{dalleCommonInterfaceAutomatic2025} and SparseMatrixColorings.jl~\cite{montoisonRevisitingSparseMatrix2025}, which are based on graph coloring~\cite{gebremedhinWhatColorYour2005}.

Companion libraries cover specialized use cases. Here, CorleoneInverse.jl implements inverse problems, currently parameter estimation with several experiments, priors declared on the annotations, and statistics of the estimates (\cref{sec:example-learning}), CorleoneOED.jl implements optimal experimental design with information-based criteria for one or several experiments (\cref{sec:example-oed}), and CorleoneGame.jl~\cite{sagerCorleoneGameLibraryBiological2026} addresses $N$-player dynamic games via generalized open-loop Nash equilibria, computed from coupled optimal control problems solved with Corleone.jl.

Problems with several stages, e.g., with different model equations, are formulated by composing one shooting layer per stage in a \texttt{Stages} container; \texttt{StageTransition} and \texttt{StageLinkage} layers implement the transition functions $F_{\mathrm{tr},i}$ of problem~\eqref{eq:problemclass}, by substitution or by coupling constraints, respectively.
Differential-algebraic and stochastic differential equations can be integrated by DifferentialEquations.jl as well, but they are not yet covered by the tests of derivatives and sparsity patterns in the test suite.

\section{Illustrative examples} \label{sec:examples}
We illustrate Corleone.jl with three variants of the Lotka--Volterra predator--prey model, covering optimal control, the training of a learned interaction term, and optimal experimental design. All three are based on tutorials in the documentation, which executes them in every build. Standalone scripts, software versions, and the hardware for the reported numbers are given in \cref{app:benchmarks}; computation times are median wall-clock times after compilation and serve only as an indication of practical cost.

\subsection{Optimal control} \label{sec:example-control}

The Lotka--Volterra fishing problem~\cite{sagerNumericalMethodsOptimal2006} asks for a fishing intensity $u(t) \in [0, 1]$ that steers the prey $x_1$ and the predators $x_2$ towards the equilibrium $(1, 1)$:
\begin{equation}
  \begin{aligned}
    \min_{u}\quad & \int_{0}^{12} \bigl(x_1(t) - 1\bigr)^2 + \bigl(x_2(t) - 1\bigr)^2 \, \mathrm{d}t\\
    \text{s.t.}\quad & \dot{x}_1 = x_1 - x_1 x_2 - 0.4\, u\, x_1, \qquad
                       \dot{x}_2 = -x_2 + x_1 x_2 - 0.2\, u\, x_2,\\
                     & x(0) = (0.5, 0.7), \qquad u(t) \in [0, 1], \qquad t \in [0, 12].
  \end{aligned}
  \label{eq:fishing}
\end{equation}
\begin{lstlisting}[language=Julia, float=tbp, caption={Lotka--Volterra fishing problem with single shooting, abridged from the Corleone.jl tutorial. \texttt{u}$_3$ is the default name of the third state in the \texttt{ODEProblem}.}, label={lst:fishing}]
using Corleone, OrdinaryDiffEqTsit5, ComponentArrays, ForwardDiff
using Optimization, OptimizationMOI, Ipopt
function lotka!(dx, x, p, t)
    dx[1] = x[1] - p[1] * x[1] * x[2] - 0.4 * p[2] * x[1]
    dx[2] = -x[2] + p[3] * x[1] * x[2] - 0.2 * p[2] * x[2]
    dx[3] = (x[1] - 1.0)^2 + (x[2] - 1.0)^2
end
prob = ODEProblem(lotka!, [0.5, 0.7, 0.0], (0.0, 12.0), [1.0, 0.0, 1.0];
                  abstol = 1e-9, reltol = 1e-9)
grid = collect(0.25:0.25:11.75)              # 48 control intervals
control = Tunable(2; kind = Val(:parameter), default = 1.0,
                  lowerbound = 0.0, upperbound = 1.0) => grid
fixed = Constant([3, 1]; kind = Val(:parameter)) => nothing
layer = ShootingLayer(prob, Tsit5(), control, fixed;
                      controls = [2], quadratures = [3])
dynprob = CorleoneDynamicOptProblem(layer, :(u₃(12.0)))
opt = OptimizationProblem(dynprob, AutoForwardDiff())
sol = solve(opt, Ipopt.Optimizer();
            hessian_approximation = "limited-memory")
\end{lstlisting}

\begin{figure}[t]
  \centering
  \begin{subfigure}[b]{0.49\linewidth}
    \centering
    \includegraphics[width=\linewidth]{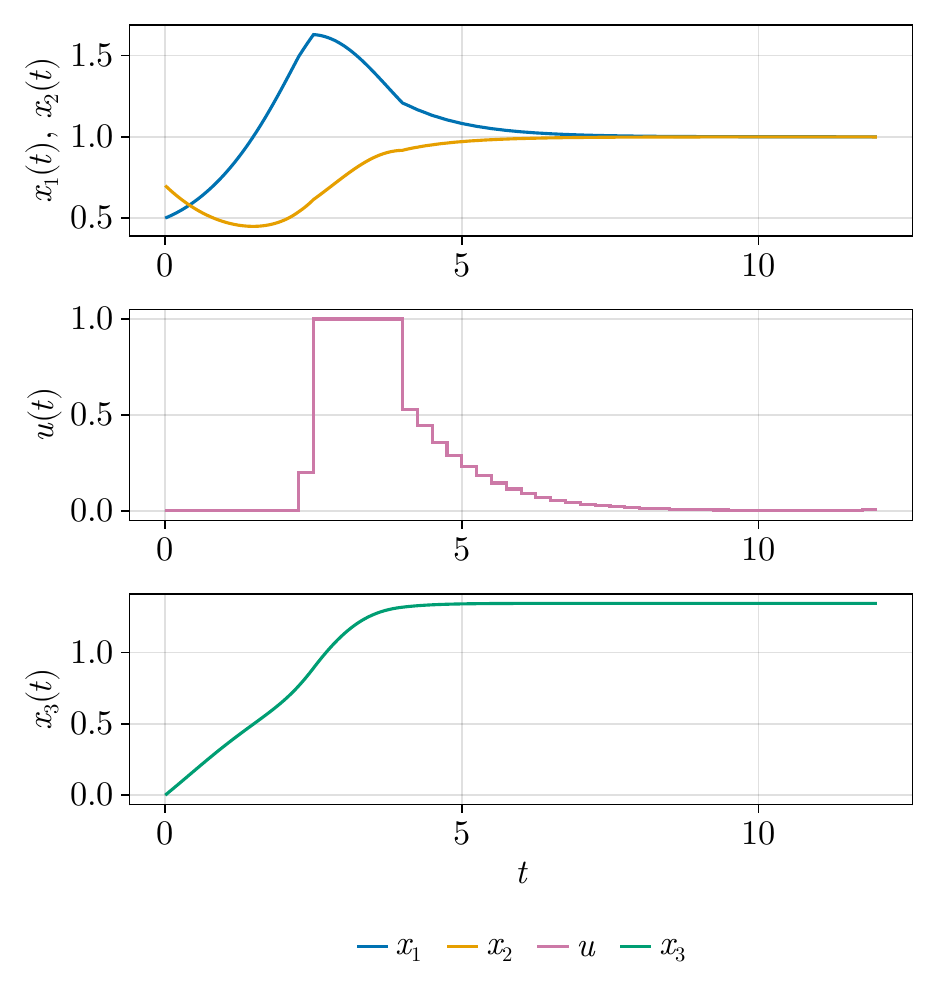}
    \caption{Single shooting, optimal solution.}
  \end{subfigure}\hfill
  \begin{subfigure}[b]{0.49\linewidth}
    \centering
    \includegraphics[width=\linewidth]{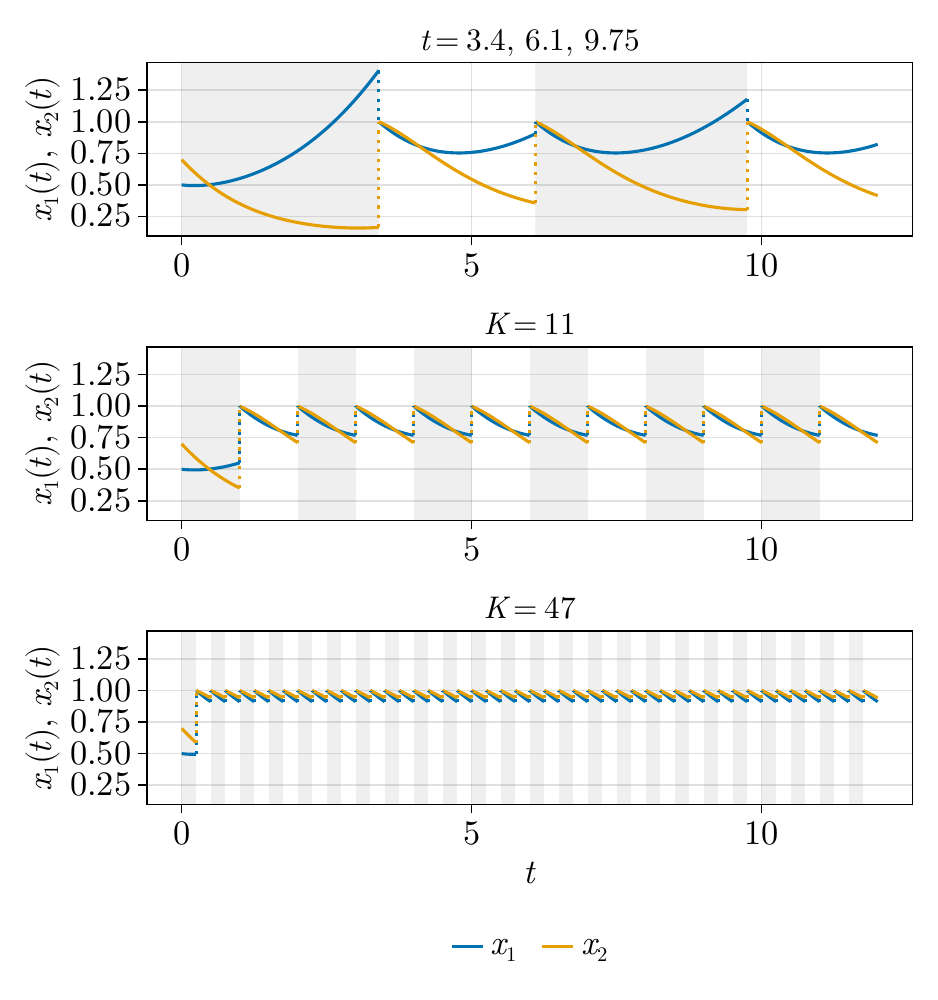}
    \caption{Multiple shooting, reference initialization.}
  \end{subfigure}
  \caption{Lotka--Volterra fishing problem. (a) Optimal solution with single shooting; from top to bottom: states $x_1, x_2$, fishing control $u$, and accumulated cost $x_3$. (b) Initial states $x_1, x_2$ of multiple shooting with the reference initialization of the lifted states; from top to bottom: shooting nodes at $t = 3.4$, $6.1$, $9.75$, and $K = 11$ and $K = 47$ equidistant nodes. Shaded areas mark the shooting blocks, dotted lines the continuity residuals at the shooting nodes.}
  \label{fig:fishing}
\end{figure}

\Cref{lst:fishing} shows the workflow. The cost integral becomes a quadrature state $x_3$; the fishing intensity is annotated as a tunable parameter and paired with a grid of 48 intervals of length $0.25$; the two interaction rates are annotated as constant; and the objective is the expression $x_3(12)$, which evaluates $x_3$ at the final time (\cref{sec:functionalities}). With single shooting, Ipopt~\cite{wachterImplementationInteriorpointFilter2006} with a limited-memory Hessian approximation converges to an objective value of $1.345633$ with Ipopt's optimality tolerance of $10^{-6}$ in 141 iterations and a runtime of 0.33\,s (\cref{fig:fishing}a). For multiple shooting, the control grid is kept, and all three states are lifted at the shooting nodes $t = 3.4$, $6.1$, and $9.75$, which need not lie on the control grid. Compared with \cref{lst:fishing}, only the lifted states and the ensemble algorithm are added:
\begin{lstlisting}[language=Julia, float=false]
lifted = (Tunable(i; kind = Val(:variable), lowerbound = 0.0,
                  upperbound = 2.0) => [3.4, 6.1, 9.75] for i in 1:3)
layer = ShootingLayer(prob, Tsit5(), control, fixed, lifted...;
    ensemblealg = EnsembleThreads(), controls = [2], quadratures = [3])
ps, st = add_parameter_segments!(layer, LuxCore.setup(rng, layer)...)
dynprob = CorleoneDynamicOptProblem(layer, :(u₃(12.0)); ps, st)
\end{lstlisting}
The last line splits the control intervals that contain the nodes $3.4$ and $6.1$, so that each shooting block owns its control values; two matching constraints keep the split values equal, which yields 59 variables and 11 constraints. Before optimization, the lifted states are initialized independently of the preceding intervals, so the trajectory jumps at the shooting nodes (\cref{fig:fishing}b); these jumps are the residuals $r_{l,k}$ that the NLP solver drives to zero. With four threads, the solver converges to an objective value of $1.345631$ with a maximum continuity residual of $6.1 \cdot 10^{-10}$ in 197 iterations and 0.42\,s. Corleone.jl also derives the sparsity patterns of the constraint Jacobian and the Lagrangian Hessian from the shooting structure. With \texttt{sparse\_constraints = true}, which requires loading DifferentiationInterface.jl and SparseMatrixColorings.jl, the Jacobian is evaluated by colored forward-mode differentiation, which reduces the runtime to 0.34\,s for the same number of iterations; exact sparse Hessians do not pay off on this small problem (\cref{tab:benchmarks}). A scaling study (\cref{tab:scaling}) shows that multiple shooting is always slower than single shooting with one thread; with a sparse Jacobian and four or more threads, it is faster for three to eleven nodes, by up to 27\,\%, although part of the difference stems from different iteration counts.

\paragraph{Shooting grid and initialization} The nodes $3.4$, $6.1$, and $9.75$ demonstrate that shooting nodes need not lie on the control grid, but they are not a good choice for convergence.
Both the shooting grid and the initialization of the lifted states have a strong influence (\cref{app:shooting}): with an exact Hessian, the iterations for these nodes drop from 51 to 23 if the lifted states are initialized at the reference
$x_1 = x_2 = 1$ of the objective instead of the initial values, and to 19 and 17 for 11 and 47 equidistant nodes (\cref{fig:fishing}b). Shooting nodes are therefore not beneficial per se: for this stable problem, single shooting converges in 14 iterations, fewer than any multiple-shooting variant.In contrast, multiple shooting is advantageous for unstable dynamics, when knowledge about the states is available for the initialization~\cite{bockplitt_1984_jul_multipleshootingalgorithm},and for the parallel integration of expensive models (\cref{tab:scaling}).

\paragraph{Integer controls} Integrality is declared in the annotation of a tunable quantity. Changing the control of \cref{lst:fishing} to 
\begin{lstlisting}[language=Julia, float=false]
Tunable(2; kind = Val(:parameter), default = 1.0, lowerbound = 0.0,
                  upperbound = 1.0, integer = true) => grid
\end{lstlisting}
restricts the fishing intensity to $u(t) \in \{0, 1\}$. This information is passed into the \texttt{OptimizationProblem} and Optimization.jl~\cite{dixitOptimizationjlUnifiedOptimization2023}, so that MINLP solvers~\cite{bonamiAlgorithmicFrameworkConvex2008} with a MathOptInterface wrapper can be used. Here, all 48 control values are marked as integer; the additional optimization variables (the lifted shooting states) are continuous. Since the problem is nonconvex and the ODE solution is not available in closed form, spatial branch-and-bound~\cite{belottiBranchingBoundsTightening2009} is not applicable. The nonlinear branch-and-bound solver Juniper~\cite{krogerJuniperOpenSourceNonlinear2018} prunes with local solutions of the NLP relaxations and therefore acts as a heuristic here; it returns an integer solution with objective value $1.352362$ after 52 branches in 83\,s, compared to $1.345633$ for the relaxed problem.

Alternatively, relaxed solutions can be rounded by the combinatorial integral approximation~\cite{sagerReformulationsAlgorithmsOptimization2009, sagerCombinatorialIntegralApproximation2011, sagerIntegerApproximationError2012}, see also~\cite{bengeaOptimalControlSwitching2005, vasudevanConsistentApproximationsOptimal2013}, and the resulting integer controls passed back as fixed controls. Sum-up rounding is currently implemented in the example scripts and is planned for the package. For the fishing problem, 33 of the 48 relaxed control values are fractional. Sum-up rounding turns them into an integer control with six switches and an objective value of $1.364083$ in 0.12\,ms, which is 1.4\,\% above the relaxed and 0.9\,\% above the MINLP solution. The approximation error decreases with the length of the control intervals~\cite{sagerIntegerApproximationError2012}, and tools such as pycombina~\cite{burger_pycombina_2020}, switching-cost-aware rounding (SCARP)~\cite{bestehornetal_2021_mp_switchingcostsshortestpath}, and rounding via matching (RVM)~\cite{bestehornkirches_2020_oo_matchingalgorithms} additionally take switching constraints or switching costs into account; a total variation penalty is an alternative~\cite{leyfferSequentialLinearInteger2022}.

\subsection{Learned dynamics} \label{sec:example-learning}

\begin{figure}[t]
  \centering
  \includegraphics[width=\linewidth]{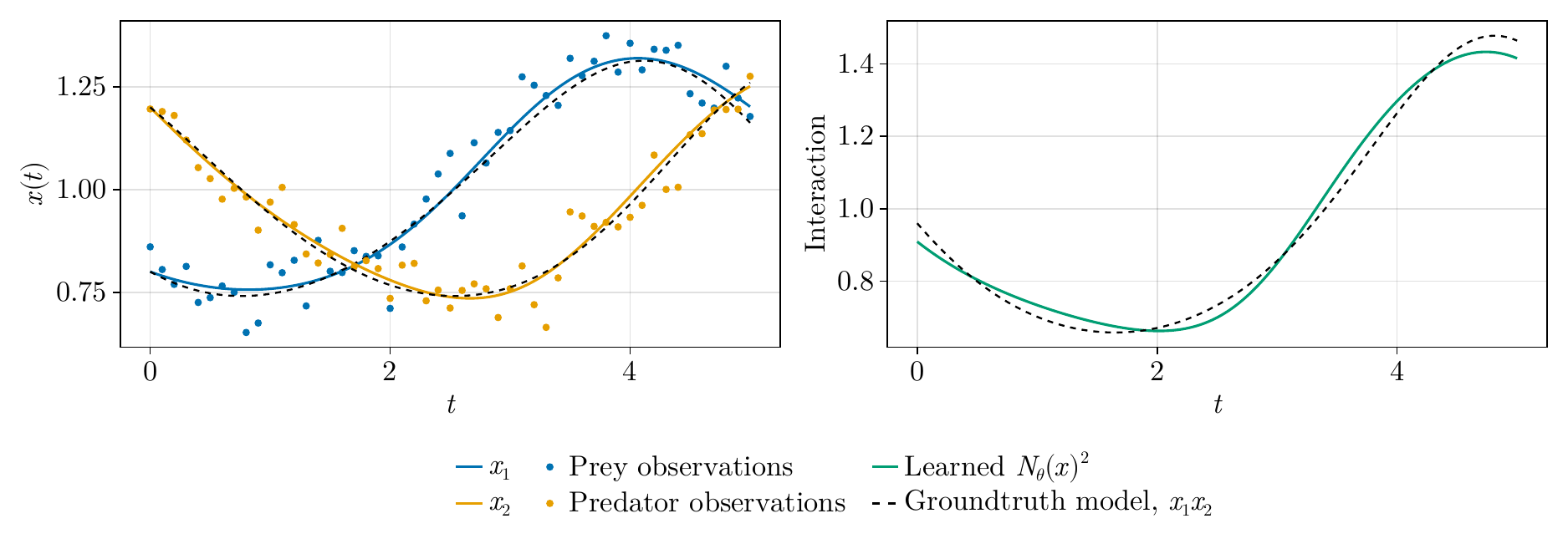}
  \caption{Trained universal differential equation. Left: states of the trained model (solid), noisy observations (dots), and the true model (dashed). Right: learned interaction $N_\theta(x)^2$ along the trajectory of the trained model and true interaction $x_1 x_2$ along the true trajectory.}
  \label{fig:ude}
\end{figure}

The second example replaces the unknown interaction term by a neural network $N_\theta$, which turns the model into a universal differential equation~\cite{rackauckasetal_2021_nov_universaldifferentialequations}:
\begin{equation}
  \begin{aligned}
    \dot{x}_1 &= x_1 - N_\theta(x)^2, & \dot{x}_2 &= N_\theta(x)^2 - x_2,\\
    x(0) &= (0.8, 1.2), & t &\in [0, 5].
  \end{aligned}
  \label{eq:ude}
\end{equation}
Squaring the network output keeps the learned interaction nonnegative. $N_\theta$ is a Lux.jl network with one hidden layer of five $\tanh$ units and 21 weights. As the network is part of the right-hand side, it is evaluated inside the \texttt{ODEProblem} that the \texttt{ShootingLayer} wraps, and its flattened weights are the parameters of this problem. In the notation of \cref{eq:ivp}, they are time-invariant parameters $c = \theta$, annotated as tunable with bounds $[-3, 3]$, so the network enters the optimization like any other parameter. Synthetic observations of both populations are generated at 51 equidistant times from the model with the true interaction $x_1 x_2$ and perturbed with Gaussian noise of standard deviation $\sigma = 0.05$. With CorleoneInverse.jl, the observations and their standard deviation form an \texttt{Experiment}, from which the weighted residuals are built automatically. The resulting nonlinear least-squares problem is solved with the Levenberg--Marquardt method~\cite{levenbergMethodSolutionCertain1944, marquardtAlgorithmLeastSquares1963} of NonlinearSolve.jl~\cite{palNonlinearSolvejlHighPerformanceRobust2024} and forward-mode Jacobians through the integrator. \Cref{lst:ude} shows the setup. Within an iteration limit of 100 iterations and 0.10\,s, at which the solver stops before convergence, the mean squared error decreases from $139.9$ to $2.53 \cdot 10^{-3}$, which is close to the noise variance $\sigma^2 = 2.5 \cdot 10^{-3}$: the network reproduces the data to the noise level (\cref{fig:ude}). As the data cover a single orbit, this does not identify the interaction on the whole state space.

\begin{lstlisting}[language=Julia, float=tbp, caption={Universal differential equation for the Lotka--Volterra model with CorleoneInverse.jl, abridged from \texttt{examples/ude.jl} of the examples repository (data generation and the scaling of the initial weights omitted). The observations refer to the states by their indices 1 and 2.}, label={lst:ude}]
network = Lux.Chain(Lux.Dense(2 => 5, tanh), Lux.Dense(5 => 1))
nn_ps, nn_st = Lux.setup(rng, network)
p0, rebuild = pack(default_vectorizer(), ComponentArray{Float64}(nn_ps))
function hybrid(x, p, t)
    out, _ = Lux.apply(network, x, rebuild(p), Lux.testmode(nn_st))
    interaction = only(out)^2
    return [x[1] - interaction, interaction - x[2]]
end
prob = ODEProblem(hybrid, [0.8, 1.2], (0.0, 5.0), p0)
weights = Tunable(collect(eachindex(p0)); kind = Val(:parameter),
                  lowerbound = -3.0, upperbound = 3.0) => nothing
layer = ShootingLayer(prob, Tsit5(), weights)
ps, st = Lux.setup(rng, layer)
experiment = Experiment(1 => (times, data[1, :], 0.05),
                        2 => (times, data[2, :], 0.05))
model = CorleoneDynamicOptProblem(layer, experiment; ps, st)
fit = solve(NonlinearLeastSquaresProblem(model),
            LevenbergMarquardt(); maxiters = 100)
\end{lstlisting}

\subsection{Optimal experimental design} \label{sec:example-oed}

\begin{figure}[t]
  \centering
  \includegraphics[width=\linewidth]{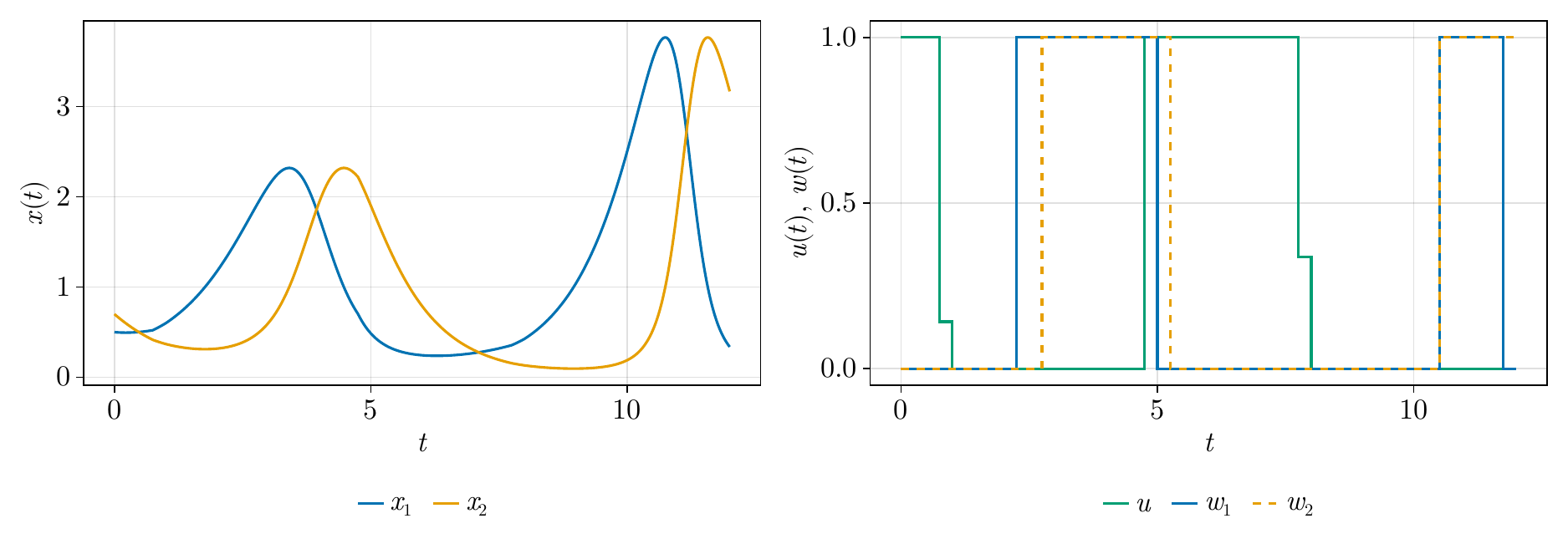}
  \caption{Optimized experimental design. Left: populations $x_1, x_2$. Right: fishing intensity $u$ and the measurement weights $w_1, w_2$ for prey and predators.}
  \label{fig:oed}
\end{figure}

The third example uses CorleoneOED.jl to design an experiment for estimating the two interaction coefficients of the fishing model in \cref{eq:fishing}. The design variables are the fishing intensity, piecewise constant on 48 intervals, and two measurement weights $w_1(t), w_2(t) \in [0, 1]$, piecewise constant on the same 48 intervals, which decide when the prey and the predator populations are observed. CorleoneOED.jl augments the \texttt{ShootingLayer} with the parameter sensitivities and the Fisher information matrix
\begin{equation}
  F = \int_{0}^{12} \sum_{m=1}^{2} w_m(t)\, J_m(t)^\top J_m(t) \, \mathrm{d}t,
  \label{eq:fim}
\end{equation}
where $J_m$ is the sensitivity of the $m$-th observation with respect to the estimated coefficients, with unit measurement variance~\cite{sager_2013_jan_samplingdecisionsoptimum}. The A-criterion $\operatorname{tr}(F^{-1})$ is minimized subject to a measurement budget of $\int_0^{12} w_m \, \mathrm{d}t \le 4$ for each population. \Cref{lst:oed} shows the setup. Ipopt returns a design with an A-criterion of $0.0370$ in 34 iterations and 0.54\,s, using the full budget of 4.0 for both populations (\cref{fig:oed}). Although the problem is relaxed, both measurement weights end up at their bounds almost everywhere, and the measurements concentrate where the populations change quickly. For comparison, observing both populations uniformly with the same budget, i.e., $w_m(t) \equiv 1/3$, and without fishing, i.e., the initial design, yields $0.161$, so the optimized fishing and measurement design reduces the criterion by a factor of 4.3.
Because the relaxed weights are at their bounds, sum-up rounding yields a binary design with practically the same A-criterion (\cref{app:oed-rounding}). The design is a local solution for the nominal parameter values.

\begin{lstlisting}[language=Julia, float=tbp, caption={Optimal experimental design for the Lotka--Volterra model with CorleoneOED.jl, abridged from the tutorial. The dynamics \texttt{lotka!} are analogous to \cref{lst:fishing}, without the cost state and with reordered parameters: \texttt{p[1]} is the fishing intensity and \texttt{p[2]}, \texttt{p[3]} are the estimated interaction coefficients.}, label={lst:oed}]
prob = ODEProblem(lotka!, [0.5, 0.7], (0.0, 12.0), [0.0, 1.0, 1.0])
fishing = Tunable(1; default = 0.0, lowerbound = 0.0,
                  upperbound = 1.0) => collect(0.25:0.25:11.75)
grid = collect(0.0:0.25:11.75)
measurements = (
    ContinuousMeasurement(:w1, grid, (x, p, t) -> x[1]),
    ContinuousMeasurement(:w2, grid, (x, p, t) -> x[2]),
)
oed = OEDLayer(prob, Tsit5(), fishing; params = [2, 3], measurements,
               controls = [1])
ps, st = LuxCore.setup(Xoshiro(1), oed)
model = CorleoneDynamicOptProblem(oed, ACriterion(); ps, st,
                                  M = [4.0, 4.0])
result = solve(OptimizationProblem(model, AutoForwardDiff()),
               Ipopt.Optimizer();
               hessian_approximation = "limited-memory")
\end{lstlisting}

\section{Impact} \label{sec:impact}

\paragraph{New research questions} Corleone.jl makes it possible to treat mechanistic, hybrid, and data-driven models within the same optimization workflow. This enables questions that require both classes of models at once, e.g., the design of experiments for universal differential equations~\cite{plateetal_2026_mar_optimalexperimentaldesign}, the computation of Nash equilibria~\cite{basarDynamicNoncooperativeGame1999} in biological dynamic games~\cite{sagerCorleoneGameLibraryBiological2026}, and, in future work, the hierarchy of inverse problems up to inverse optimal control and inverse dynamic games~\cite{abazarInverseOptimalControl2020, molloyInverseOptimalControl2022, bangaGeneralizedInverseOptimal2025, bangaHierarchicalFrameworkInverse2026}.

\paragraph{Improvement of existing research} Users choose integrators, derivative backends, and optimizers independently, and they can switch between single and multiple shooting by changing the annotations rather than the model. The companion libraries show that specialized methods require comparatively little additional code, since they reuse the shooting layer, the trajectory interface, and the optimization bridge: the source code of CorleoneOED.jl comprises about 1,100, that of CorleoneGame.jl about 1,000, and that of CorleoneInverse.jl about 700 non-blank, non-comment lines of Julia, compared to about 5,500 for the core of Corleone.jl (v2.1.0).

\paragraph{Daily practice} The predecessor package DynamicOED.jl~\cite{martensenDynamicOEDjlJuliaPackage2024}, which was used for the computations in~\cite{plateetal_2026_mar_optimalexperimentaldesign}, has been merged into Corleone.jl and CorleoneOED.jl.
Corleone.jl is used in approximately one dozen ongoing master's and PhD projects and in the courses ``Optimization Methods for Machine Learning'' and ``Modeling, Simulation and Optimization'' at Otto von Guericke University Magdeburg. The executable tutorials serve as starting points for students and new users.

\paragraph{Spread of use and commercial use} Corleone.jl has been used in prototype studies in the research departments of BASF for optimal experimental design and of Volkswagen for hybrid modeling. Reliable user numbers are not yet available.

\section{Conclusions} \label{sec:conclusion}

Corleone.jl represents the piecewise simulation of a dynamical system as a Lux.jl layer and thereby acts as a transcription layer between dynamic optimization problems and (mixed-integer) nonlinear programming solvers. Mechanistic and hybrid models, integrators, derivative backends, and solvers can be combined freely, and companion libraries add optimal experimental design, inverse problems, and dynamic games with little additional code. Future work includes inverse optimal control and inverse dynamic games, dedicated support for mixed-integer optimal control, starting with sum-up rounding in the package, and an extension of the maintained derivative coverage to differential-algebraic and stochastic differential equations.

\section*{Acknowledgements}
\label{sec:acknowledgements}

\begin{sloppypar}
This work was supported by 
the Deutsche Forschungsgemeinschaft (DFG, German Research Foundation) via the priority program 2331 ``Machine Learning in Chemical Engineering'' under grant SA 2016/3-1;
the European Regional Development Fund under the European Union's Horizon Europe Research and Innovation Program [grant numbers ZS/2023/12/182064 (intelAlgen), ZS/2023/12/182063 (timingMatters), and ZS/2023/12/182075 (Center for Dynamic Systems)];
and BASF SE. 
We thank Hergen Schultze, Chryssa Kappatou, Robert Lee, and Kai Exner of BASF's Data Analytics group for discussions and industrial use cases.
The funding sources had no involvement in the writing of the article. All software resulting from this work, i.e., Corleone.jl and its companion libraries CorleoneOED.jl, CorleoneInverse.jl, and CorleoneGame.jl, is open source and freely available under the MIT license.
\end{sloppypar}

We thank the Julia SciML community, whose ecosystem of differential equation solvers, sensitivity analysis, automatic differentiation, and optimization interfaces Corleone.jl builds upon.

\section*{\contributionsheading}

\begin{sloppypar}
\textbf{Carl Julius Martensen}: Conceptualization, Methodology, Software (Corleone.jl, CorleoneOED.jl, CorleoneInverse.jl, CorleoneGame.jl), Validation, Formal analysis, Writing -- original draft, Writing -- review \& editing.\\
\textbf{Christoph Plate}: Conceptualization, Methodology, Software (Corleone.jl, CorleoneOED.jl, CorleoneInverse.jl, CorleoneGame.jl), Validation, Formal analysis, Writing -- original draft, Writing -- review \& editing.\\
\textbf{Christian Kirches}: Writing -- review \& editing.\\
\textbf{Andreas Potschka}: Writing -- review \& editing.\\
\textbf{Christopher Vincent Rackauckas}: Writing -- review \& editing.\\
\textbf{Sebastian Sager}: Conceptualization, Methodology, Software (CorleoneOED.jl, CorleoneGame.jl), Writing -- original draft, Writing -- review \& editing, Supervision, Funding acquisition.
\end{sloppypar}

\section*{\competingheading}
\label{sec:competing}

This work was partly funded by BASF SE, in whose research department Corleone.jl has been used in prototype studies (\cref{sec:impact}). The authors declare that they have no other known competing financial interests or personal relationships that could have appeared to influence the work reported in this paper.

\section*{\aiheading}

During the preparation of this work the authors used Claude Code 2.1.285 with the model Claude Opus 5.5 in order to revise code, search the literature, and draft and revise parts of the text. After using this tool, the authors reviewed and edited the content as needed and take full responsibility for the content of \thisarticle.

\appendix
\section{Benchmark of the illustrative examples} \label{app:benchmarks}

\Cref{tab:benchmarks} reports the results and computation times of the examples of \cref{sec:examples}, obtained with the standalone scripts in the \texttt{examples} folder of \url{https://github.com/mathopt/Corleone-Paper-Examples}~\cite{martensenetal_2026_corleonepaperexamples}. The scripts state the problems as in \cref{lst:fishing,lst:ude,lst:oed} and use only registered releases of the General Julia registry; the repository fixes the environment in its \texttt{Project.toml} and \texttt{Manifest.toml}, so that \texttt{Pkg.instantiate()} restores the versions of \cref{tab:dependencies}, and its \texttt{results} folder contains the timings, logs, and figures of the run reported here. Ipopt runs with a tolerance of $10^{-6}$. The times are wall-clock times of a complete solve, measured with BenchmarkTools.jl after a warm-up run, so compilation and problem construction are excluded.

\begin{table}[H]
  \centering
  \small
  \setlength{\tabcolsep}{4pt}
  \begin{tabular}{@{}llrrll@{}}
    \toprule
    Example & Solver & Var. & Iter. & Objective & Time [s] \\
    \midrule
    Single shooting              & Ipopt    &  48 & 141 & $1.345633$           & $0.33$ \\
    Integer controls             & Juniper  &  48 & --  & $1.352362$           & $83.4$ \\
    \quad sum-up rounding        & --       &  48 & --  & $1.364083$           & $1.2 \cdot 10^{-4}$ \\
    Multiple shooting            & Ipopt    &  59 & 197 & $1.345631$           & $0.42$ \\
    \quad sparse $J$             & Ipopt    &  59 & 197 & $1.345631$           & $0.34$ \\
    \quad sparse $J$ and $H$     & Ipopt    &  59 &  51 & $1.345633$           & $0.60$ \\
    UDE                          & LM       &  21 & 100 & $2.53 \cdot 10^{-3}$ & $0.10$ \\
    OED                          & Ipopt    & 144 &  34 & $3.70 \cdot 10^{-2}$ & $0.54$ \\
    \quad sum-up rounding        & --       &  96 & --  & $3.70 \cdot 10^{-2}$ & $2.2 \cdot 10^{-4}$ \\
    \bottomrule
  \end{tabular}
  \caption{Results and computation times of the examples of \cref{sec:examples} with the software of \cref{tab:dependencies} and four threads. Time: median wall-clock time of 10 samples. Juniper uses Ipopt for the relaxations; LM: Levenberg--Marquardt. Var.: number of optimization variables; Iter.: solver iterations; Objective: $x_3(12)$ for the fishing problems, mean squared error for the UDE, and A-criterion $\operatorname{tr}(F^{-1})$ for the OED. Multiple shooting reaches a maximum continuity residual of $6.1 \cdot 10^{-10}$. Its sparse variants evaluate the constraint Jacobian $J$ ($11 \times 59$, 154 nonzeros) and, in addition, the exact Lagrangian Hessian $H$ ($59 \times 59$, 889 nonzeros) with the sparsity patterns derived by Corleone.jl and colored forward-mode differentiation; the latter replaces the L-BFGS approximation of the other Ipopt runs and needs fewer but more expensive iterations. They reach maximum continuity residuals of $7.0 \cdot 10^{-10}$ and $2.4 \cdot 10^{-7}$. Juniper does not report iterations and needs 52 branches. The UDE fit stops at its iteration limit of 100. The sum-up rounding rows round the relaxed solutions of the single-shooting fishing problem (48 control values) and of the OED (96 measurement weights); they are timed together with the evaluation of the rounded solution, without the relaxed solve.}
  \label{tab:benchmarks}
\end{table}

\paragraph{Rounding of the optimal experimental design} \label{app:oed-rounding} The script \texttt{oed\_cia.jl} rounds the relaxed measurement weights of \cref{sec:example-oed} by sum-up rounding~\cite{sagerReformulationsAlgorithmsOptimization2009, sagerCombinatorialIntegralApproximation2011}. Here, the relaxation is tight: the relaxed weights are binary up to $1.5 \cdot 10^{-6}$, with 16 of the 48 intervals selected per population. The rounded design uses the budgets of exactly 4, and its A-criterion of $0.0370$ differs from the relaxed one by less than $10^{-9}$. In contrast, the relaxed fishing control has 33 fractional values on 48 intervals, and its rounding increases the objective by 1.4\,\% (\cref{sec:example-control}).

{\sloppy\paragraph{Scaling of multiple shooting} \Cref{tab:scaling} reports the solve times of the fishing problem for different numbers of shooting nodes and threads for a threaded ensemble. The script \texttt{scaling.sh} runs \texttt{fishing\_multiple\_shooting\_scaling.jl} for each number of threads. The problem is that of \cref{sec:example-control} with the same 48 control intervals, but with $K$ equidistant shooting nodes, at which all three states are lifted; $K = 0$ is single shooting. For $K = 3$, the nodes lie at $t = 3$, $6$, and $9$ instead of $3.4$, $6.1$, and $9.75$, so the iteration counts differ from \cref{tab:benchmarks}. Each configuration is solved with a dense and a sparse constraint Jacobian and the limited-memory Hessian approximation of Ipopt, and timed with five samples after a warm-up run. All runs converge to an objective value of $1.345631$ ($1.345633$ for single shooting), and the iteration counts do not depend on the number of threads. They differ between the dense and the sparse Jacobian for some $K$ because the limited-memory approximation reacts to small numerical differences in the derivatives. In addition, the time per iteration depends on the number of segments: the threaded ensemble integrates each of the $K + 1$ segments in its own task, so $K = 1$ keeps at most two threads busy, whereas $K = 3$ distributes the integration over four, which makes an iteration with four or more threads 34--44\,\% cheaper. With one thread, multiple shooting is always slower than single shooting; with four or more threads and a sparse Jacobian, it is faster for three to eleven nodes, whereas with a dense Jacobian only $K = 3$, and $K = 5$ with eight threads, are faster. As the iteration counts differ between configurations, the per-iteration times $t/\text{iter.}$ isolate the effect of the threads: with eight threads and $K = 5$, the sparse variant needs 1.9\,ms per iteration, compared to 2.4\,ms for single shooting. In our tests, which are not part of the repository, distributed worker processes (\texttt{EnsembleDistributed()}) returned the same solutions, but every evaluation transfers the problem and the dual numbers of the forward-mode differentiation to the workers; for this small model, this overhead is not worth it, and they are therefore not included. They are meant for expensive blocks of stages, e.g., of large or stiff models.\par}

\begin{table}[H]
  \centering
  \small
  \setlength{\tabcolsep}{4pt}
  \begin{tabular}{@{}rrlcccc@{}}
    \toprule
    $K$ & Var. & & 1 thread & 2 threads & 4 threads & 8 threads \\
    \midrule
     0 &  48 & time [s] & 0.33 & 0.34 & 0.33 & 0.34 \\
           &     & \footnotesize iter. & \footnotesize 141 & \footnotesize 141 & \footnotesize 141 & \footnotesize 141 \\[2pt]
     1 &  51 & time [s] & 0.74 / 0.52 & 0.54 / \textbf{0.34} & 0.53 / 0.34 & 0.72 / 0.46 \\
           &     & \footnotesize iter. & \footnotesize 177 / 134 & \footnotesize 177 / 134 & \footnotesize 177 / 134 & \footnotesize 177 / 134 \\[2pt]
     3 &  57 & time [s] & 0.48 / 0.54 & 0.34 / 0.38 & \textbf{0.24} / \textbf{0.27} & \textbf{0.30} / \textbf{0.32} \\
           &     & \footnotesize iter. & \footnotesize 120 / 167 & \footnotesize 120 / 167 & \footnotesize 120 / 167 & \footnotesize 120 / 167 \\[2pt]
     5 &  63 & time [s] & 0.59 / 0.48 & 0.43 / \textbf{0.32} & 0.36 / \textbf{0.26} & \textbf{0.31} / \textbf{0.25} \\
           &     & \footnotesize iter. & \footnotesize 129 / 128 & \footnotesize 129 / 128 & \footnotesize 129 / 128 & \footnotesize 129 / 128 \\[2pt]
    11 &  81 & time [s] & 0.76 / 0.51 & 0.58 / \textbf{0.34} & 0.43 / \textbf{0.27} & 0.43 / \textbf{0.26} \\
           &     & \footnotesize iter. & \footnotesize 136 / 120 & \footnotesize 136 / 120 & \footnotesize 136 / 120 & \footnotesize 136 / 120 \\[2pt]
    23 & 117 & time [s] & 1.06 / 1.11 & 0.80 / 0.86 & 0.60 / 0.65 & 0.54 / 0.58 \\
           &     & \footnotesize iter. & \footnotesize 96 / 175 & \footnotesize 96 / 175 & \footnotesize 96 / 175 & \footnotesize 96 / 175 \\[2pt]
    47 & 189 & time [s] & 2.12 / 1.39 & 1.65 / 0.94 & 1.26 / 0.77 & 1.12 / 0.61 \\
           &     & \footnotesize iter. & \footnotesize 82 / 82 & \footnotesize 82 / 82 & \footnotesize 82 / 82 & \footnotesize 82 / 82 \\
    \bottomrule
  \end{tabular}
  \caption{Scaling of multiple shooting for the fishing problem with $K$ equidistant shooting nodes and a threaded ensemble. Entries: median wall-clock time of a solve (time) and Ipopt iterations (iter.), each with dense / sparse constraint Jacobian; single shooting ($K = 0$) has no constraints. Var.: number of optimization variables. Bold: faster than single shooting with the same number of threads. The fastest configuration ($K = 3$, dense, 4 threads) is 28\,\% faster than single shooting with the same number of threads.}
  \label{tab:scaling}
\end{table}

\paragraph{Shooting grid and initialization} \label{app:shooting} \Cref{tab:shooting} compares the shooting nodes $3.4$, $6.1$, and $9.75$ of \cref{sec:example-control} with single shooting and with $K = 11$ and $K = 47$ equidistant nodes, i.e., at every fourth and at every point of the control grid. Without a default value, the lifted states are initialized with the initial values $x(0) = (0.5, 0.7, 0)$ of the problem; alternatively, \mbox{\texttt{default = 1.0}} in their annotations initializes $x_1$ and $x_2$ at the reference $x_1 = x_2 = 1$ of the objective (\cref{fig:fishing}b). Since the iteration counts of the limited-memory approximation depend on small numerical differences in the derivatives, this comparison uses the exact sparse Lagrangian Hessian; all variants converge to the objective value $1.34563$. Both the grid and the initialization have a strong influence on convergence: the reference initialization halves the number of iterations on the coarse grids, and the finest grid needs the fewest iterations. The script \texttt{fishing\_shooting\_grids.jl} of the examples repository reproduces \cref{tab:shooting}.

\begin{table}[H]
  \centering
  \small
  \begin{tabular}{@{}lcccc@{}}
    \toprule
    & Single shooting & $t = 3.4, 6.1, 9.75$ & $K = 11$ & $K = 47$ \\
    \midrule
    Variables & 48 & 59 & 81 & 189 \\
    Iterations, initial values & 14 & 51 & 47 & 20 \\
    Iterations, reference & -- & 23 & 19 & 17 \\
    \bottomrule
  \end{tabular}
  \caption{Fishing problem with different shooting grids and initializations of the lifted states: number of optimization variables and Ipopt iterations with the exact sparse Lagrangian Hessian. The lifted states are initialized with the initial values of the problem or with the reference $x_1 = x_2 = 1$ of the objective.}
  \label{tab:shooting}
\end{table}

\paragraph{Software} \Cref{tab:dependencies} lists all direct dependencies of the scripts and selected indirect dependencies, as resolved in \texttt{Manifest.toml} (403 packages in total, including the plotting packages for the figures). The standard libraries are those of Julia~1.12.7.

\begin{table}[!htbp]
  \centering
  \footnotesize
  \begin{tabular}{@{}ll>{\raggedright\arraybackslash}p{5.2cm}@{}}
    \toprule
    Package & Version & Role \\
    \midrule
    \multicolumn{3}{@{}l}{\emph{Direct dependencies}} \\
    Corleone.jl                  & 2.1.0    & shooting layers, dynamic optimization problems \\
    CorleoneOED.jl               & 2.1.0    & optimal experimental design (\cref{sec:example-oed}) \\
    CorleoneInverse.jl           & 0.1.0    & experiments, least-squares problems (\cref{sec:example-learning}) \\
    OrdinaryDiffEqTsit5.jl       & 2.1.5    & Tsit5 integrator \\
    Lux.jl                       & 1.32.1   & neural network of the UDE \\
    LuxCore.jl                   & 1.5.5    & layer interface \\
    ComponentArrays.jl           & 0.15.50  & default parameter vectorization \\
    ForwardDiff.jl               & 1.5.0    & forward-mode automatic differentiation \\
    Optimization.jl              & 5.9.1    & optimization interface \\
    OptimizationMOI.jl           & 1.6.1    & MathOptInterface bridge for Ipopt and Juniper \\
    Ipopt.jl                     & 1.16.0   & interior-point NLP solver \\
    Juniper.jl                   & 0.9.5    & nonlinear branch and bound (integer controls) \\
    NonlinearSolve.jl            & 4.32.0   & Levenberg--Marquardt method (UDE) \\
    DifferentiationInterface.jl  & 0.7.22   & AD backend interface, sparse derivatives \\
    SparseMatrixColorings.jl     & 0.4.29   & coloring of sparse Jacobians and Hessians \\
    SymbolicIndexingInterface.jl & 0.3.56   & symbolic access to trajectories \\
    BenchmarkTools.jl            & 1.8.0    & timings \\
    CairoMakie.jl, LaTeXStrings.jl & 0.15.17, 1.4.1 & figures \\
    InteractiveUtils, Printf, Random & stdlib & system information, output, random numbers \\
    \midrule
    \multicolumn{3}{@{}l}{\emph{Selected indirect dependencies}} \\
    SciMLBase.jl                 & 3.58.0   & problem and solution interface \\
    DiffEqBase.jl                & 7.21.6   & differential equation solver core \\
    OrdinaryDiffEqCore.jl        & 4.19.1   & ODE integrator core \\
    OptimizationBase.jl          & 5.10.1   & derivatives for Optimization.jl \\
    NonlinearSolveFirstOrder.jl  & 2.10.0   & first-order nonlinear least-squares methods \\
    MathOptInterface.jl          & 1.54.0   & solver abstraction \\
    ADTypes.jl                   & 1.24.0   & AD backend selection \\
    RecursiveArrayTools.jl       & 4.5.6    & trajectory containers \\
    LuxLib.jl                    & 1.15.12  & neural network kernels \\
    Distributions.jl             & 0.25.134 & noise models and priors \\
    Ipopt\_jll                   & 300.1400.1902+0 & Ipopt binary \\
    MUMPS\_seq\_jll              & 500.900.100+0   & linear solver of Ipopt \\
    OpenBLAS32\_jll              & 0.3.34+1 & BLAS and LAPACK of Ipopt \\
    \bottomrule
  \end{tabular}
  \caption{Software used for \cref{tab:benchmarks}: Julia packages and versions from the \texttt{Manifest.toml} of the examples repository.}
  \label{tab:dependencies}
\end{table}

\paragraph{System} Output of \texttt{InteractiveUtils.versioninfo()} on the machine used for \cref{tab:benchmarks}:

\begin{small}
\begin{verbatim}
Julia Version 1.12.7
Commit 6d172b025e4 (2026-08-15 08:05 UTC)
Build Info:
  Official https://julialang.org release
Platform Info:
  OS: Linux (x86_64-linux-gnu)
  CPU: 16 × AMD Ryzen 7 8700F 8-Core Processor
  WORD_SIZE: 64
  LLVM: libLLVM-18.1.7 (ORCJIT, znver3)
  GC: Built with stock GC
Threads: 4 default, 1 interactive, 4 GC (on 16 virtual cores)
\end{verbatim}
\end{small}

\section{Software availability} \label{app:availability}

The results of this article were obtained with the following software versions and resources.

\begin{description}\setlength{\itemsep}{0pt}
\item[Code version] Corleone.jl v2.1.0, CorleoneOED.jl 2.1.0, and CorleoneInverse.jl 0.1.0.
\item[Repository] \url{https://github.com/SciML/Corleone.jl}, release \url{https://github.com/SciML/Corleone.jl/releases/tag/v2.1.0}; the companion libraries are tagged \texttt{CorleoneOED-v2.1.0} and \texttt{CorleoneInverse-v0.1.0} in the same repository.
\item[Examples and results] \url{https://github.com/mathopt/Corleone-Paper-Examples}, archived as v1.0.3 at \url{https://doi.org/10.5281/zenodo.24580692}~\cite{martensenetal_2026_corleonepaperexamples}.
\item[License] MIT.
\item[Languages and tools] Julia; SciML (OrdinaryDiffEq.jl, SciMLSensitivity.jl), Lux.jl, Optimization.jl, MathOptInterface.jl, and Ipopt; version control with git.
\item[Requirements] Julia $\geq$ 1.10 (CorleoneGame.jl: $\geq$ 1.11; results of this article: Julia 1.12.7); dependencies are listed in the \texttt{Project.toml} files.
\item[Documentation] \url{https://sciml.github.io/Corleone.jl/dev/}.
\item[Contact] \texttt{carl.martensen@ovgu.de}.
\end{description}

\bibliographystyle{unsrtnat}
\bibliography{references}

\begin{thebibliography}{92}
\providecommand{\natexlab}[1]{#1}
\providecommand{\url}[1]{\texttt{#1}}
\expandafter\ifx\csname urlstyle\endcsname\relax
  \providecommand{\doi}[1]{doi: #1}\else
  \providecommand{\doi}{doi: \begingroup \urlstyle{rm}\Url}\fi

\bibitem[Srinivasan et~al.(2003)Srinivasan, Palanki, and
  Bonvin]{srinivasanDynamicOptimizationBatch2003}
B.~Srinivasan, S.~Palanki, and D.~Bonvin.
\newblock Dynamic optimization of batch processes: {I}. {Characterization} of
  the nominal solution.
\newblock \emph{Computers \& Chemical Engineering}, 27\penalty0 (1):\penalty0
  1--26, 2003.
\newblock \doi{10.1016/S0098-1354(02)00116-3}.

\bibitem[Hrovat et~al.(2012)Hrovat, Di~Cairano, Tseng, and
  Kolmanovsky]{hrovatDevelopmentModelPredictive2012}
Davor Hrovat, Stefano Di~Cairano, H.~Eric Tseng, and Ilya~V. Kolmanovsky.
\newblock The development of model predictive control in automotive industry: A
  survey.
\newblock In \emph{2012 IEEE International Conference on Control Applications},
  pages 295--302, 2012.
\newblock \doi{10.1109/CCA.2012.6402735}.

\bibitem[Moore(2018)]{mooreHowMathematicallyOptimize2018}
Helen Moore.
\newblock How to mathematically optimize drug regimens using optimal control.
\newblock \emph{Journal of Pharmacokinetics and Pharmacodynamics}, 45\penalty0
  (1):\penalty0 127--137, 2018.
\newblock \doi{10.1007/s10928-018-9568-y}.

\bibitem[Raue et~al.(2013)Raue, Schilling, Bachmann, Matteson, Schelke,
  Kaschek, Hug, Kreutz, Harms, Theis, Klingm{\"u}ller, and
  Timmer]{raueLessonsLearnedQuantitative2013}
Andreas Raue, Marcel Schilling, Julie Bachmann, Andrew Matteson, Max Schelke,
  Daniel Kaschek, Sabine Hug, Clemens Kreutz, Brian~D. Harms, Fabian~J. Theis,
  Ursula Klingm{\"u}ller, and Jens Timmer.
\newblock Lessons learned from quantitative dynamical modeling in systems
  biology.
\newblock \emph{PLoS ONE}, 8\penalty0 (9):\penalty0 e74335, 2013.
\newblock \doi{10.1371/journal.pone.0074335}.

\bibitem[Wensing et~al.(2024)Wensing, Posa, Hu, Escande, Mansard, and
  Del~Prete]{wensingOptimizationbasedControlDynamic2024}
Patrick~M. Wensing, Michael Posa, Yue Hu, Adrien Escande, Nicolas Mansard, and
  Andrea Del~Prete.
\newblock Optimization-based control for dynamic legged robots.
\newblock \emph{IEEE Transactions on Robotics}, 40:\penalty0 43--63, 2024.
\newblock \doi{10.1109/TRO.2023.3324580}.

\bibitem[Hu et~al.(2021)Hu, Shan, Guerrero, Ioinovici, Chan, and
  Rodriguez]{huModelPredictiveControl2021}
Jiefeng Hu, Yinghao Shan, Josep~M. Guerrero, Adrian Ioinovici, Ka~Wing Chan,
  and Jose Rodriguez.
\newblock Model predictive control of microgrids -- an overview.
\newblock \emph{Renewable and Sustainable Energy Reviews}, 136:\penalty0
  110422, 2021.
\newblock \doi{10.1016/j.rser.2020.110422}.

\bibitem[Betts(2010)]{bettsPracticalMethodsOptimal2010}
John~T. Betts.
\newblock \emph{Practical Methods for Optimal Control and Estimation Using
  Nonlinear Programming}.
\newblock Advances in Design and Control. SIAM, Philadelphia, 2 edition, 2010.
\newblock \doi{10.1137/1.9780898718577}.

\bibitem[Biegler(2010)]{bieglerNonlinearProgrammingConcepts2010}
Lorenz~T. Biegler.
\newblock \emph{Nonlinear Programming: Concepts, Algorithms, and Applications
  to Chemical Processes}.
\newblock MOS-SIAM Series on Optimization. SIAM, Philadelphia, 2010.
\newblock \doi{10.1137/1.9780898719383}.

\bibitem[Sass and Mitsos(2025)]{sassObscuredTerminologyHidden2025}
Susanne Sass and Alexander Mitsos.
\newblock Obscured by terminology: {Hidden} parallels in direct methods for
  open-loop optimal control.
\newblock \emph{Journal of Process Control}, 154:\penalty0 103513, October
  2025.
\newblock \doi{10.1016/j.jprocont.2025.103513}.

\bibitem[Bock and Plitt(1984)]{bockplitt_1984_jul_multipleshootingalgorithm}
H.G. Bock and K.J. Plitt.
\newblock A {Multiple} {Shooting} {Algorithm} for {Direct} {Solution} of
  {Optimal} {Control} {Problems}.
\newblock \emph{IFAC Proceedings Volumes}, 17\penalty0 (2):\penalty0
  1603--1608, July 1984.
\newblock ISSN 1474-6670.
\newblock \doi{10.1016/S1474-6670(17)61205-9}.

\bibitem[Leineweber et~al.(2003)Leineweber, Bauer, Bock, and
  Schl{\"o}der]{leineweberEfficientMultipleShooting2003}
Daniel~B. Leineweber, Irene Bauer, Hans~Georg Bock, and Johannes~P.
  Schl{\"o}der.
\newblock An efficient multiple shooting based reduced {SQP} strategy for
  large-scale dynamic process optimization. {Part} 1: Theoretical aspects.
\newblock \emph{Computers \& Chemical Engineering}, 27\penalty0 (2):\penalty0
  157--166, 2003.
\newblock \doi{10.1016/S0098-1354(02)00158-8}.

\bibitem[Rao et~al.(1998)Rao, Wright, and
  Rawlings]{raoApplicationInteriorpointMethods1998}
Christopher~V. Rao, Stephen~J. Wright, and James~B. Rawlings.
\newblock Application of interior-point methods to model predictive control.
\newblock \emph{Journal of Optimization Theory and Applications}, 99\penalty0
  (3):\penalty0 723--757, 1998.
\newblock \doi{10.1023/A:1021711402723}.

\bibitem[Frison and Diehl(2020)]{frisonHPIPMHighperformanceQuadratic2020}
Gianluca Frison and Moritz Diehl.
\newblock {HPIPM}: a high-performance quadratic programming framework for model
  predictive control.
\newblock \emph{IFAC-PapersOnLine}, 53\penalty0 (2):\penalty0 6563--6569, 2020.
\newblock \doi{10.1016/j.ifacol.2020.12.073}.

\bibitem[Janka et~al.(2016)Janka, Kirches, Sager, and
  Wächter]{jankaetal_2016_dec_sr1bfgssqp}
Dennis Janka, Christian Kirches, Sebastian Sager, and Andreas Wächter.
\newblock An {SR1}/{BFGS} {SQP} algorithm for nonconvex nonlinear programs with
  block-diagonal {Hessian} matrix.
\newblock \emph{Mathematical Programming Computation}, 8\penalty0 (4):\penalty0
  435--459, December 2016.
\newblock ISSN 1867-2949, 1867-2957.
\newblock \doi{10.1007/s12532-016-0101-2}.

\bibitem[Wittmann and Sager(2025)]{wittmannsager___blocksqp2exploiting}
Reinhold Wittmann and Sebastian Sager.
\newblock {blockSQP} 2: exploiting structure to improve performance.
\newblock Preprint, Optimization Online, October 2025.
\newblock URL
  \url{https://optimization-online.org/2025/10/blocksqp-2-exploiting-structure-to-improve-performance/}.

\bibitem[Bürger et~al.(2020)Bürger, Zeile, Hahn, Altmann-Dieses, Sager, and
  Diehl]{burger_pycombina_2020}
Adrian Bürger, Clemens Zeile, Mirko Hahn, Angelika Altmann-Dieses, Sebastian
  Sager, and Moritz Diehl.
\newblock pycombina: {An} {Open}-{Source} {Tool} for {Solving} {Combinatorial}
  {Approximation} {Problems} {Arising} in {Mixed}-{Integer} {Optimal}
  {Control}.
\newblock \emph{IFAC-PapersOnLine}, 53\penalty0 (2):\penalty0 6502--6508, 2020.
\newblock ISSN 24058963.
\newblock \doi{10.1016/j.ifacol.2020.12.1799}.

\bibitem[Bauer et~al.(2000)Bauer, Bock, K{\"o}rkel, and
  Schl{\"o}der]{bauerNumericalMethodsOptimum2000}
Irene Bauer, Hans~Georg Bock, Stefan K{\"o}rkel, and Johannes~P. Schl{\"o}der.
\newblock Numerical methods for optimum experimental design in {DAE} systems.
\newblock \emph{Journal of Computational and Applied Mathematics}, 120\penalty0
  (1--2):\penalty0 1--25, 2000.
\newblock \doi{10.1016/S0377-0427(00)00300-9}.

\bibitem[Franceschini and
  Macchietto(2008)]{franceschiniModelbasedDesignExperiments2008}
Gaia Franceschini and Sandro Macchietto.
\newblock Model-based design of experiments for parameter precision: State of
  the art.
\newblock \emph{Chemical Engineering Science}, 63\penalty0 (19):\penalty0
  4846--4872, 2008.
\newblock \doi{10.1016/j.ces.2007.11.034}.

\bibitem[Sager(2013)]{sager_2013_jan_samplingdecisionsoptimum}
Sebastian Sager.
\newblock Sampling {Decisions} in {Optimum} {Experimental} {Design} in the
  {Light} of {Pontryagin}'s {Maximum} {Principle}.
\newblock \emph{SIAM Journal on Control and Optimization}, 51\penalty0
  (4):\penalty0 3181--3207, January 2013.
\newblock ISSN 0363-0129, 1095-7138.
\newblock \doi{10.1137/110835098}.

\bibitem[Wang and Dowling(2022)]{wangPyomoDOEOpensourcePackage2022}
Jialu Wang and Alexander~W. Dowling.
\newblock Pyomo.{DOE}: {An} open‐source package for model‐based design of
  experiments in {Python}.
\newblock \emph{AIChE Journal}, 68\penalty0 (12):\penalty0 e17813, December
  2022.
\newblock \doi{10.1002/aic.17813}.

\bibitem[Martensen et~al.(2024)Martensen, Plate, and
  Sager]{martensenDynamicOEDjlJuliaPackage2024}
Carl~Julius Martensen, Christoph Plate, and Sebastian Sager.
\newblock {DynamicOED}.jl: {A} {Julia} package for solving optimum experimental
  design problems.
\newblock \emph{Journal of Open Source Software}, 9\penalty0 (98):\penalty0
  6605, June 2024.
\newblock ISSN 2475-9066.
\newblock \doi{10.21105/joss.06605}.

\bibitem[Karniadakis et~al.(2021)Karniadakis, Kevrekidis, Lu, Perdikaris, Wang,
  and Yang]{karniadakisPhysicsinformedMachineLearning2021}
George~Em Karniadakis, Ioannis~G. Kevrekidis, Lu~Lu, Paris Perdikaris, Sifan
  Wang, and Liu Yang.
\newblock Physics-informed machine learning.
\newblock \emph{Nature Reviews Physics}, 3\penalty0 (6):\penalty0 422--440,
  2021.
\newblock \doi{10.1038/s42254-021-00314-5}.

\bibitem[Rackauckas et~al.(2021)Rackauckas, Ma, Martensen, Warner, Zubov,
  Supekar, Skinner, Ramadhan, and
  Edelman]{rackauckasetal_2021_nov_universaldifferentialequations}
Christopher Rackauckas, Yingbo Ma, Carl~Julius Martensen, Collin Warner, Kirill
  Zubov, Rohit Supekar, Dominic Skinner, Ali Ramadhan, and Alan Edelman.
\newblock Universal {Differential} {Equations} for {Scientific} {Machine}
  {Learning}, November 2021.
\newblock arXiv:2001.04385 [cs, math, q-bio, stat].

\bibitem[Plate et~al.(2026)Plate, Martensen, and
  Sager]{plateetal_2026_mar_optimalexperimentaldesign}
Christoph Plate, Carl~Julius Martensen, and Sebastian Sager.
\newblock Optimal {Experimental} {Design} for {Universal} {Differential}
  {Equations}.
\newblock \emph{IEEE Transactions on Automatic Control}, 71\penalty0
  (3):\penalty0 1521--1536, March 2026.
\newblock ISSN 0018-9286, 1558-2523, 2334-3303.
\newblock \doi{10.1109/TAC.2025.3609533}.

\bibitem[Andersson et~al.(2019)Andersson, Gillis, Horn, Rawlings, and
  Diehl]{anderssonCasADiSoftwareFramework2019}
Joel A.~E. Andersson, Joris Gillis, Greg Horn, James~B. Rawlings, and Moritz
  Diehl.
\newblock {CasADi}: a software framework for nonlinear optimization and optimal
  control.
\newblock \emph{Mathematical Programming Computation}, 11\penalty0
  (1):\penalty0 1--36, March 2019.
\newblock \doi{10.1007/s12532-018-0139-4}.

\bibitem[Nicholson et~al.(2018)Nicholson, Siirola, Watson, Zavala, and
  Biegler]{nicholsonPyomodaeModelingAutomatic2018}
Bethany Nicholson, John~D. Siirola, Jean-Paul Watson, Victor~M. Zavala, and
  Lorenz~T. Biegler.
\newblock pyomo.dae: a modeling and automatic discretization framework for
  optimization with differential and algebraic equations.
\newblock \emph{Mathematical Programming Computation}, 10\penalty0
  (2):\penalty0 187--223, June 2018.
\newblock \doi{10.1007/s12532-017-0127-0}.

\bibitem[Pulsipher et~al.(2022)Pulsipher, Zhang, Hongisto, and
  Zavala]{pulsipherUnifyingModelingAbstraction2022}
Joshua~L. Pulsipher, Weiqi Zhang, Tyler~J. Hongisto, and Victor~M. Zavala.
\newblock A unifying modeling abstraction for infinite-dimensional
  optimization.
\newblock \emph{Computers \& Chemical Engineering}, 156:\penalty0 107567,
  January 2022.
\newblock \doi{10.1016/j.compchemeng.2021.107567}.

\bibitem[Caillau et~al.(2026)Caillau, Cots, Gergaud, Martinon, and
  Sed]{caillauOptimalControljlJuliaPackage2026}
Jean-Baptiste Caillau, Olivier Cots, Joseph Gergaud, Pierre Martinon, and
  Sophia Sed.
\newblock {OptimalControl}.jl: a {Julia} package to model and solve optimal
  control problems with {ODE}'s.
\newblock Zenodo, April 2026.
\newblock Version 2.0.3.

\bibitem[Ma et~al.(2021)Ma, Gowda, Anantharaman, Laughman, Shah, and
  Rackauckas]{maModelingToolkitComposableGraph2021}
Yingbo Ma, Shashi Gowda, Ranjan Anantharaman, Chris Laughman, Viral Shah, and
  Chris Rackauckas.
\newblock {ModelingToolkit}: a composable graph transformation system for
  equation-based modeling, March 2021.
\newblock arXiv: 2103.05244.

\bibitem[Salzmann et~al.(2024)Salzmann, Arrizabalaga, Andersson, Pavone, and
  Ryll]{salzmannLearningCasADiDatadriven2024}
Tim Salzmann, Jon Arrizabalaga, Joel Andersson, Marco Pavone, and Markus Ryll.
\newblock Learning for {CasADi}: {Data}-driven models in numerical
  optimization.
\newblock In \emph{Proceedings of the 6th annual learning for dynamics \&
  control conference}, volume 242 of \emph{Proceedings of machine learning
  research}, pages 541--553, 2024.
\newblock URL \url{https://proceedings.mlr.press/v242/salzmann24a.html}.

\bibitem[Lahr et~al.(2024)Lahr, Näf, Wabersich, Frey, Siehl, Carron, Diehl,
  and Zeilinger]{lahrL4acadosLearningbasedModels2024}
Amon Lahr, Joshua Näf, Kim~P. Wabersich, Jonathan Frey, Pascal Siehl, Andrea
  Carron, Moritz Diehl, and Melanie~N. Zeilinger.
\newblock L4acados: {Learning}-based models for acados, applied to {Gaussian}
  process-based predictive control, November 2024.
\newblock arXiv: 2411.19258.

\bibitem[Verschueren et~al.(2022)Verschueren, Frison, Kouzoupis, Frey,
  Duijkeren, Zanelli, Novoselnik, Albin, Quirynen, and
  Diehl]{verschuerenAcadosModularOpensource2022}
Robin Verschueren, Gianluca Frison, Dimitris Kouzoupis, Jonathan Frey,
  Niels~van Duijkeren, Andrea Zanelli, Branimir Novoselnik, Thivaharan Albin,
  Rien Quirynen, and Moritz Diehl.
\newblock acados—a modular open-source framework for fast embedded optimal
  control.
\newblock \emph{Mathematical Programming Computation}, 14\penalty0
  (1):\penalty0 147--183, March 2022.
\newblock \doi{10.1007/s12532-021-00208-8}.

\bibitem[Fr{\"o}hlich et~al.(2021)Fr{\"o}hlich, Weindl, Sch{\"a}lte, Pathirana,
  Paszkowski, Lines, Stapor, and
  Hasenauer]{frohlichAMICIHighperformanceSensitivity2021}
Fabian Fr{\"o}hlich, Daniel Weindl, Yannik Sch{\"a}lte, Dilan Pathirana,
  {\L}ukasz Paszkowski, Glenn~Terje Lines, Paul Stapor, and Jan Hasenauer.
\newblock {AMICI}: high-performance sensitivity analysis for large ordinary
  differential equation models.
\newblock \emph{Bioinformatics}, 37\penalty0 (20):\penalty0 3676--3677, 2021.
\newblock \doi{10.1093/bioinformatics/btab227}.

\bibitem[Sch{\"a}lte et~al.(2023)Sch{\"a}lte, Fr{\"o}hlich, Jost, Vanhoefer,
  Pathirana, Stapor, Lakrisenko, Wang, Raim{\'u}ndez, Merkt, Schmiester,
  St{\"a}dter, Grein, Dudkin, Doresic, Weindl, and
  Hasenauer]{schaltePyPESTOModularScalable2023}
Yannik Sch{\"a}lte, Fabian Fr{\"o}hlich, Paul~J. Jost, Jakob Vanhoefer, Dilan
  Pathirana, Paul Stapor, Polina Lakrisenko, Dantong Wang, Elba Raim{\'u}ndez,
  Simon Merkt, Leonard Schmiester, Philipp St{\"a}dter, Stephan Grein, Erika
  Dudkin, Domagoj Doresic, Daniel Weindl, and Jan Hasenauer.
\newblock {pyPESTO}: a modular and scalable tool for parameter estimation for
  dynamic models.
\newblock \emph{Bioinformatics}, 39\penalty0 (11):\penalty0 btad711, 2023.
\newblock \doi{10.1093/bioinformatics/btad711}.

\bibitem[Schmiester et~al.(2021)Schmiester, Sch{\"a}lte, Bergmann, Camba,
  Dudkin, Egert, Fr{\"o}hlich, Fuhrmann, Hauber, Kemmer, Lakrisenko, Loos,
  Merkt, M{\"u}ller, Pathirana, Raim{\'u}ndez, Refisch, Rosenblatt, Stapor,
  St{\"a}dter, Wang, Wieland, Banga, Timmer, Villaverde, Sahle, Kreutz,
  Hasenauer, and Weindl]{schmiesterPEtabInteroperableSpecification2021}
Leonard Schmiester, Yannik Sch{\"a}lte, Frank~T. Bergmann, Tacio Camba, Erika
  Dudkin, Janine Egert, Fabian Fr{\"o}hlich, Lara Fuhrmann, Adrian~L. Hauber,
  Svenja Kemmer, Polina Lakrisenko, Carolin Loos, Simon Merkt, Wolfgang
  M{\"u}ller, Dilan Pathirana, Elba Raim{\'u}ndez, Lukas Refisch, Marcus
  Rosenblatt, Paul~L. Stapor, Philipp St{\"a}dter, Dantong Wang, Franz-Georg
  Wieland, Julio~R. Banga, Jens Timmer, Alejandro~F. Villaverde, Sven Sahle,
  Clemens Kreutz, Jan Hasenauer, and Daniel Weindl.
\newblock {PEtab}---interoperable specification of parameter estimation
  problems in systems biology.
\newblock \emph{PLOS Computational Biology}, 17\penalty0 (1):\penalty0
  e1008646, 2021.
\newblock \doi{10.1371/journal.pcbi.1008646}.

\bibitem[Balsa-Canto et~al.(2016)Balsa-Canto, Henriques, G{\'a}bor, and
  Banga]{balsacantoAMIGO2ToolboxDynamic2016}
Eva Balsa-Canto, David Henriques, Attila G{\'a}bor, and Julio~R. Banga.
\newblock {AMIGO2}, a toolbox for dynamic modeling, optimization and control in
  systems biology.
\newblock \emph{Bioinformatics}, 32\penalty0 (21):\penalty0 3357--3359, 2016.
\newblock \doi{10.1093/bioinformatics/btw411}.

\bibitem[Bonnans et~al.(2012)Bonnans, Martinon, and
  Gr{\'e}lard]{bonnansBocopCollectionExamples2012}
Fr{\'e}d{\'e}ric~J. Bonnans, Pierre Martinon, and Vincent Gr{\'e}lard.
\newblock {BOCOP} -- {A} collection of examples.
\newblock Research Report RR-8053, Inria, 2012.
\newblock {HAL}: hal-00726992.

\bibitem[Fiedler et~al.(2023)Fiedler, Karg, L{\"u}ken, Brandner, Heinlein,
  Brabender, and Lucia]{fiedlerDompcFAIRNonlinear2023}
Felix Fiedler, Benjamin Karg, Lukas L{\"u}ken, Dean Brandner, Moritz Heinlein,
  Felix Brabender, and Sergio Lucia.
\newblock do-mpc: Towards {FAIR} nonlinear and robust model predictive control.
\newblock \emph{Control Engineering Practice}, 140:\penalty0 105676, 2023.
\newblock \doi{10.1016/j.conengprac.2023.105676}.

\bibitem[Falck et~al.(2021)Falck, Gray, Ponnapalli, and
  Wright]{falckDymosPythonPackage2021}
Robert Falck, Justin~S. Gray, Kaushik Ponnapalli, and Ted Wright.
\newblock dymos: A {Python} package for optimal control of multidisciplinary
  systems.
\newblock \emph{Journal of Open Source Software}, 6\penalty0 (59):\penalty0
  2809, 2021.
\newblock \doi{10.21105/joss.02809}.

\bibitem[Caspari et~al.(2019)Caspari, Bremen, Faust, Jung, Kappatou, Sass,
  Vaupel, Hannemann-Tam{\'a}s, Mhamdi, and
  Mitsos]{caspariDyOSFrameworkOptimization2019}
Adrian Caspari, Andreas~M. Bremen, Johannes M.~M. Faust, Florian Jung,
  Chrysoula~D. Kappatou, Susanne Sass, Yannic Vaupel, Ralf Hannemann-Tam{\'a}s,
  Adel Mhamdi, and Alexander Mitsos.
\newblock {DyOS} -- a framework for optimization of large-scale differential
  algebraic equation systems.
\newblock In \emph{29th European Symposium on Computer Aided Process
  Engineering}, volume~46 of \emph{Computer Aided Chemical Engineering}, pages
  619--624. Elsevier, 2019.
\newblock \doi{10.1016/B978-0-12-818634-3.50104-1}.

\bibitem[Beal et~al.(2018)Beal, Hill, Martin, and
  Hedengren]{bealGEKKOOptimizationSuite2018}
Logan D.~R. Beal, Daniel~C. Hill, R.~Abraham Martin, and John~D. Hedengren.
\newblock {GEKKO} optimization suite.
\newblock \emph{Processes}, 6\penalty0 (8):\penalty0 106, 2018.
\newblock \doi{10.3390/pr6080106}.

\bibitem[Patterson and Rao(2014)]{pattersonGPOPSIIMATLABSoftware2014}
Michael~A. Patterson and Anil~V. Rao.
\newblock {GPOPS-II}: A {MATLAB} software for solving multiple-phase optimal
  control problems using hp-adaptive {Gaussian} quadrature collocation methods
  and sparse nonlinear programming.
\newblock \emph{ACM Transactions on Mathematical Software}, 41\penalty0
  (1):\penalty0 1--37, 2014.
\newblock \doi{10.1145/2558904}.

\bibitem[Agamawi and Rao(2020)]{agamawiCGPOPSOptimalControl2020}
Yunus~M. Agamawi and Anil~V. Rao.
\newblock {CGPOPS}: A {C++} software for solving multiple-phase optimal control
  problems using adaptive {Gaussian} quadrature collocation and sparse
  nonlinear programming.
\newblock \emph{ACM Transactions on Mathematical Software}, 46\penalty0
  (3):\penalty0 1--38, 2020.
\newblock \doi{10.1145/3390463}.

\bibitem[Pohlodek et~al.(2025)Pohlodek, Morabito, Schlauch, Zometa, and
  Findeisen]{pohlodekFlexibleDevelopmentEvaluation2025}
Johannes Pohlodek, Bruno Morabito, Christian Schlauch, Pablo Zometa, and Rolf
  Findeisen.
\newblock Flexible development and evaluation of machine-learning-supported
  optimal control and estimation methods via {HILO-MPC}.
\newblock \emph{International Journal of Robust and Nonlinear Control},
  35\penalty0 (7):\penalty0 2835--2859, 2025.
\newblock \doi{10.1002/rnc.7275}.

\bibitem[Nie et~al.(2018)Nie, Faqir, and Kerrigan]{nieICLOCS2TryThis2018}
Yuanbo Nie, Omar Faqir, and Eric~C. Kerrigan.
\newblock {ICLOCS2}: Try this optimal control problem solver before you try the
  rest.
\newblock In \emph{2018 UKACC 12th International Conference on Control
  (CONTROL)}, page 336, Sheffield, UK, 2018.
\newblock \doi{10.1109/CONTROL.2018.8516795}.

\bibitem[Hoffmann et~al.(2024)Hoffmann, Kirches, Potschka, Sager, and
  Wirsching]{hoffmannMUSCODIIUserManual2024}
Christian Hoffmann, Christian Kirches, Andreas Potschka, Sebastian Sager, and
  Leonard Wirsching.
\newblock \emph{{MUSCOD-II} User Manual, Release 6.0}.
\newblock Interdisciplinary Center for Scientific Computing (IWR), University
  of Heidelberg, 2024.
\newblock Based on the previous version by M. Diehl, D. B. Leineweber, and A.
  A. S. Sch{\"a}fer.

\bibitem[Drgo{\v{n}}a et~al.(2022)Drgo{\v{n}}a, Ki{\v{s}}, Tuor, Vrabie, and
  Klau{\v{c}}o]{drgonaDifferentiablePredictiveControl2022}
J{\'a}n Drgo{\v{n}}a, Karol Ki{\v{s}}, Aaron Tuor, Draguna Vrabie, and Martin
  Klau{\v{c}}o.
\newblock Differentiable predictive control: Deep learning alternative to
  explicit model predictive control for unknown nonlinear systems.
\newblock \emph{Journal of Process Control}, 116:\penalty0 80--92, 2022.
\newblock \doi{10.1016/j.jprocont.2022.06.001}.

\bibitem[Persson et~al.(2025)Persson, Fr{\"o}hlich, Grein, Loman, Ognissanti,
  Hasselgren, Hasenauer, and Cvijovic]{perssonPEtabjlAdvancingEfficiency2025}
Sebastian Persson, Fabian Fr{\"o}hlich, Stephan Grein, Torkel Loman, Damiano
  Ognissanti, Viktor Hasselgren, Jan Hasenauer, and Marija Cvijovic.
\newblock {PEtab.jl}: advancing the efficiency and utility of dynamic
  modelling.
\newblock \emph{Bioinformatics}, 41\penalty0 (9):\penalty0 btaf497, 2025.
\newblock \doi{10.1093/bioinformatics/btaf497}.

\bibitem[Becerra(2010)]{becerraSolvingComplexOptimal2010}
Victor~M. Becerra.
\newblock Solving complex optimal control problems at no cost with {PSOPT}.
\newblock In \emph{2010 IEEE International Symposium on Computer-Aided Control
  System Design}, pages 1391--1396, 2010.
\newblock \doi{10.1109/CACSD.2010.5612676}.

\bibitem[Chowdhary et~al.(2024)Chowdhary, Ahmed, and
  Attia]{chowdharyPyOEDExtensibleSuite2024}
Abhijit Chowdhary, Shady~E. Ahmed, and Ahmed Attia.
\newblock {PyOED}: An extensible suite for data assimilation and
  model-constrained optimal design of experiments.
\newblock \emph{ACM Transactions on Mathematical Software}, 50\penalty0
  (2):\penalty0 1--22, 2024.
\newblock \doi{10.1145/3653071}.

\bibitem[Hart et~al.(2011)Hart, Watson, and
  Woodruff]{hartPyomoModelingSolving2011}
William~E. Hart, Jean-Paul Watson, and David~L. Woodruff.
\newblock Pyomo: Modeling and solving mathematical programs in {Python}.
\newblock \emph{Mathematical Programming Computation}, 3\penalty0 (3):\penalty0
  219--260, 2011.
\newblock \doi{10.1007/s12532-011-0026-8}.

\bibitem[K{\"o}rkel(2002)]{koerkelNumerischeMethodenOptimale2002}
Stefan K{\"o}rkel.
\newblock \emph{{Numerische Methoden f{\"u}r Optimale Versuchsplanungsprobleme
  bei nichtlinearen DAE-Modellen}}.
\newblock PhD thesis, Universit{\"a}t Heidelberg, Heidelberg, 2002.

\bibitem[Bezanson et~al.(2017)Bezanson, Edelman, Karpinski, and
  Shah]{bezansonJuliaFreshApproach2017}
Jeff Bezanson, Alan Edelman, Stefan Karpinski, and Viral~B. Shah.
\newblock Julia: {A} {Fresh} {Approach} to {Numerical} {Computing}.
\newblock \emph{SIAM Review}, 59\penalty0 (1):\penalty0 65--98, January 2017.
\newblock ISSN 0036-1445.
\newblock \doi{10.1137/141000671}.

\bibitem[Pal(2023)]{palLuxExplicitParameterization2023}
Avik Pal.
\newblock Lux: {Explicit} parameterization of deep neural networks in {Julia}.
\newblock Zenodo, April 2023.

\bibitem[Rackauckas and
  Nie(2017)]{rackauckasDifferentialEquationsjlPerformantFeaturerich2017}
Christopher Rackauckas and Qing Nie.
\newblock {DifferentialEquations}.jl – a performant and feature-rich
  ecosystem for solving differential equations in {Julia}.
\newblock \emph{Journal of Open Research Software}, 5\penalty0 (1), 2017.
\newblock \doi{10.5334/jors.151}.

\bibitem[Dixit and
  Rackauckas(2023)]{dixitOptimizationjlUnifiedOptimization2023}
Vaibhav~Kumar Dixit and Christopher Rackauckas.
\newblock Optimization.jl: a unified optimization package.
\newblock Zenodo, March 2023.
\newblock Version 3.12.1.

\bibitem[Blochwitz et~al.(2011)Blochwitz, Otter, Arnold, Bausch, Clau{\ss},
  Elmqvist, Junghanns, Mauss, Monteiro, Neidhold,
  et~al.]{blochwitzFunctionalMockupInterface2011}
Torsten Blochwitz, Martin Otter, M.~Arnold, C.~Bausch, C.~Clau{\ss},
  H.~Elmqvist, A.~Junghanns, J.~Mauss, M.~Monteiro, T.~Neidhold, et~al.
\newblock The {Functional} {Mockup} {Interface} for tool independent exchange
  of simulation models.
\newblock In \emph{Proceedings of the 8th International {Modelica} Conference},
  pages 105--114, Dresden, Germany, 2011. Link{\"o}ping University Electronic
  Press.
\newblock \doi{10.3384/ecp11063105}.

\bibitem[Psichogios and Ungar(1992)]{psichogiosHybridNeuralNetworkfirst1992}
Dimitris~C. Psichogios and Lyle~H. Ungar.
\newblock A hybrid neural network--first principles approach to process
  modeling.
\newblock \emph{AIChE Journal}, 38\penalty0 (10):\penalty0 1499--1511, 1992.
\newblock \doi{10.1002/aic.690381003}.

\bibitem[Schweidtmann et~al.(2024)Schweidtmann, Zhang, and von
  Stosch]{schweidtmannReviewPerspectiveHybrid2024}
Artur~M. Schweidtmann, Dongda Zhang, and Moritz von Stosch.
\newblock A review and perspective on hybrid modeling methodologies.
\newblock \emph{Digital Chemical Engineering}, 10:\penalty0 100136, 2024.
\newblock \doi{10.1016/j.dche.2023.100136}.

\bibitem[Chen et~al.(2018)Chen, Rubanova, Bettencourt, and
  Duvenaud]{chenNeuralOrdinaryDifferential2018}
Ricky T.~Q. Chen, Yulia Rubanova, Jesse Bettencourt, and David Duvenaud.
\newblock Neural ordinary differential equations.
\newblock In \emph{Advances in Neural Information Processing Systems 31
  (NeurIPS 2018)}, 2018.

\bibitem[Kidger(2021)]{kidgerNeuralDifferentialEquations2021}
Patrick Kidger.
\newblock \emph{On Neural Differential Equations}.
\newblock PhD thesis, University of Oxford, 2021.
\newblock arXiv:2202.02435.

\bibitem[Massaroli et~al.(2021)Massaroli, Poli, Sonoda, Suzuki, Park,
  Yamashita, and Asama]{massaroliDifferentiableMultipleShooting2021}
Stefano Massaroli, Michael Poli, Sho Sonoda, Taiji Suzuki, Jinkyoo Park,
  Atsushi Yamashita, and Hajime Asama.
\newblock Differentiable multiple shooting layers.
\newblock In \emph{Advances in Neural Information Processing Systems 34
  (NeurIPS 2021)}, 2021.
\newblock URL
  \url{https://proceedings.neurips.cc/paper/2021/hash/89b9c689a57b82e59074c6ba09aa394d-Abstract.html}.

\bibitem[Turan and J{\"a}schke(2022)]{turanMultipleShootingTraining2022}
Evren~Mert Turan and Johannes J{\"a}schke.
\newblock Multiple shooting for training neural differential equations on time
  series.
\newblock \emph{IEEE Control Systems Letters}, 6:\penalty0 1897--1902, 2022.
\newblock \doi{10.1109/LCSYS.2021.3135835}.

\bibitem[Revels et~al.(2016)Revels, Lubin, and
  Papamarkou]{revelsForwardModeAutomaticDifferentiation2016}
Jarrett Revels, Miles Lubin, and Theodore Papamarkou.
\newblock Forward-mode automatic differentiation in {Julia}, 2016.

\bibitem[Innes(2018)]{innesDonUnrollAdjoint2018}
Michael Innes.
\newblock Don't unroll adjoint: Differentiating {SSA}-form programs, 2018.

\bibitem[Sapienza et~al.(2026)Sapienza, Bolibar, Sch{\"a}fer, Groenke, Pal,
  Boussange, Heimbach, Hooker, P{\'e}rez, Persson, and
  Rackauckas]{sapienza_2026_diffeqad}
Facundo Sapienza, Jordi Bolibar, Frank Sch{\"a}fer, Brian Groenke, Avik Pal,
  Victor Boussange, Patrick Heimbach, Giles Hooker, Fernando P{\'e}rez,
  Per-Olof Persson, and Christopher Rackauckas.
\newblock Differentiable programming for differential equations: A review.
\newblock \emph{SIAM Review}, 68\penalty0 (4):\penalty0 749--819, 2026.
\newblock \doi{10.1137/24M1669062}.
\newblock URL \url{https://doi.org/10.1137/24M1669062}.

\bibitem[Dalle and Hill(2025)]{dalleCommonInterfaceAutomatic2025}
Guillaume Dalle and Adrian Hill.
\newblock A common interface for automatic differentiation, 2025.

\bibitem[Montoison et~al.(2025)Montoison, Dalle, and
  Gebremedhin]{montoisonRevisitingSparseMatrix2025}
Alexis Montoison, Guillaume Dalle, and Assefaw Gebremedhin.
\newblock Revisiting sparse matrix coloring and bicoloring, 2025.

\bibitem[Gebremedhin et~al.(2005)Gebremedhin, Manne, and
  Pothen]{gebremedhinWhatColorYour2005}
Assefaw~Hadish Gebremedhin, Fredrik Manne, and Alex Pothen.
\newblock What color is your {Jacobian}? {Graph} coloring for computing
  derivatives.
\newblock \emph{SIAM Review}, 47\penalty0 (4):\penalty0 629--705, 2005.
\newblock \doi{10.1137/S0036144504444711}.

\bibitem[Sager et~al.(2026)Sager, Plate, Martensen, and
  Ohl]{sagerCorleoneGameLibraryBiological2026}
Sebastian Sager, Christoph Plate, Carl~Julius Martensen, and Frank~W. Ohl.
\newblock {CorleoneGame}: A library of biological dynamic games with a
  companion {Julia} solver, October 2026.

\bibitem[Sager et~al.(2006)Sager, Bock, Diehl, Reinelt, and
  Schl{\"o}der]{sagerNumericalMethodsOptimal2006}
Sebastian Sager, Hans~Georg Bock, Moritz Diehl, Gerhard Reinelt, and
  Johannes~P. Schl{\"o}der.
\newblock Numerical methods for optimal control with binary control functions
  applied to a {Lotka--Volterra} type fishing problem.
\newblock In \emph{Recent Advances in Optimization}, volume 563 of
  \emph{Lecture Notes in Economics and Mathematical Systems}, pages 269--289.
  Springer, 2006.
\newblock \doi{10.1007/3-540-28258-0_17}.

\bibitem[W{\"a}chter and
  Biegler(2006)]{wachterImplementationInteriorpointFilter2006}
Andreas W{\"a}chter and Lorenz~T. Biegler.
\newblock On the implementation of an interior-point filter line-search
  algorithm for large-scale nonlinear programming.
\newblock \emph{Mathematical Programming}, 106\penalty0 (1):\penalty0 25--57,
  2006.
\newblock \doi{10.1007/s10107-004-0559-y}.

\bibitem[Bonami et~al.(2008)Bonami, Biegler, Conn, Cornu{\'e}jols, Grossmann,
  Laird, Lee, Lodi, Margot, Sawaya, and
  W{\"a}chter]{bonamiAlgorithmicFrameworkConvex2008}
Pierre Bonami, Lorenz~T. Biegler, Andrew~R. Conn, G{\'e}rard Cornu{\'e}jols,
  Ignacio~E. Grossmann, Carl~D. Laird, Jon Lee, Andrea Lodi, Fran{\c{c}}ois
  Margot, Nicolas Sawaya, and Andreas W{\"a}chter.
\newblock An algorithmic framework for convex mixed integer nonlinear programs.
\newblock \emph{Discrete Optimization}, 5\penalty0 (2):\penalty0 186--204,
  2008.
\newblock \doi{10.1016/j.disopt.2006.10.011}.

\bibitem[Belotti et~al.(2009)Belotti, Lee, Liberti, Margot, and
  W{\"a}chter]{belottiBranchingBoundsTightening2009}
Pietro Belotti, Jon Lee, Leo Liberti, Fran{\c{c}}ois Margot, and Andreas
  W{\"a}chter.
\newblock Branching and bounds tightening techniques for non-convex {MINLP}.
\newblock \emph{Optimization Methods and Software}, 24\penalty0
  (4--5):\penalty0 597--634, 2009.
\newblock \doi{10.1080/10556780903087124}.

\bibitem[Kr{\"o}ger et~al.(2018)Kr{\"o}ger, Coffrin, Hijazi, and
  Nagarajan]{krogerJuniperOpenSourceNonlinear2018}
Ole Kr{\"o}ger, Carleton Coffrin, Hassan Hijazi, and Harsha Nagarajan.
\newblock Juniper: An open-source nonlinear branch-and-bound solver in {Julia}.
\newblock In \emph{Integration of Constraint Programming, Artificial
  Intelligence, and Operations Research}, volume 10848 of \emph{Lecture Notes
  in Computer Science}, pages 377--386. Springer International Publishing,
  2018.
\newblock \doi{10.1007/978-3-319-93031-2_27}.

\bibitem[Sager(2009)]{sagerReformulationsAlgorithmsOptimization2009}
Sebastian Sager.
\newblock Reformulations and algorithms for the optimization of switching
  decisions in nonlinear optimal control.
\newblock \emph{Journal of Process Control}, 19\penalty0 (8):\penalty0
  1238--1247, 2009.
\newblock \doi{10.1016/j.jprocont.2009.03.008}.

\bibitem[Sager et~al.(2011)Sager, Jung, and
  Kirches]{sagerCombinatorialIntegralApproximation2011}
Sebastian Sager, Michael Jung, and Christian Kirches.
\newblock Combinatorial integral approximation.
\newblock \emph{Mathematical Methods of Operations Research}, 73\penalty0
  (3):\penalty0 363--380, 2011.
\newblock \doi{10.1007/s00186-011-0355-4}.

\bibitem[Sager et~al.(2012)Sager, Bock, and
  Diehl]{sagerIntegerApproximationError2012}
Sebastian Sager, Hans~Georg Bock, and Moritz Diehl.
\newblock The integer approximation error in mixed-integer optimal control.
\newblock \emph{Mathematical Programming}, 133\penalty0 (1--2):\penalty0 1--23,
  2012.
\newblock \doi{10.1007/s10107-010-0405-3}.

\bibitem[Bengea and DeCarlo(2005)]{bengeaOptimalControlSwitching2005}
Sorin~C. Bengea and Raymond~A. DeCarlo.
\newblock Optimal control of switching systems.
\newblock \emph{Automatica}, 41\penalty0 (1):\penalty0 11--27, 2005.
\newblock \doi{10.1016/j.automatica.2004.08.003}.

\bibitem[Vasudevan et~al.(2013)Vasudevan, Gonzalez, Bajcsy, and
  Sastry]{vasudevanConsistentApproximationsOptimal2013}
Ramanarayan Vasudevan, Humberto Gonzalez, Ruzena Bajcsy, and S.~Shankar Sastry.
\newblock Consistent approximations for the optimal control of constrained
  switched systems---{Part} 1: A conceptual algorithm.
\newblock \emph{SIAM Journal on Control and Optimization}, 51\penalty0
  (6):\penalty0 4463--4483, 2013.
\newblock \doi{10.1137/120901490}.

\bibitem[Bestehorn et~al.(2021)Bestehorn, Hansknecht, Kirches, and
  Manns]{bestehornetal_2021_mp_switchingcostsshortestpath}
Felix Bestehorn, Christoph Hansknecht, Christian Kirches, and Paul Manns.
\newblock Mixed-integer optimal control problems with switching costs: a
  shortest path approach.
\newblock \emph{Mathematical Programming}, 188\penalty0 (2):\penalty0 621--652,
  2021.
\newblock \doi{10.1007/s10107-020-01581-3}.

\bibitem[Bestehorn and
  Kirches(2020)]{bestehornkirches_2020_oo_matchingalgorithms}
Felix Bestehorn and Christian Kirches.
\newblock Matching algorithms and complexity results for constrained
  mixed-integer optimal control with switching costs.
\newblock Optimization Online preprint, 2020.
\newblock URL \url{https://optimization-online.org/?p=16748}.

\bibitem[Leyffer and Manns(2022)]{leyfferSequentialLinearInteger2022}
Sven Leyffer and Paul Manns.
\newblock Sequential linear integer programming for integer optimal control
  with total variation regularization.
\newblock \emph{ESAIM: Control, Optimisation and Calculus of Variations},
  28:\penalty0 66, 2022.
\newblock \doi{10.1051/cocv/2022059}.

\bibitem[Levenberg(1944)]{levenbergMethodSolutionCertain1944}
Kenneth Levenberg.
\newblock A method for the solution of certain non-linear problems in least
  squares.
\newblock \emph{Quarterly of Applied Mathematics}, 2\penalty0 (2):\penalty0
  164--168, 1944.
\newblock \doi{10.1090/qam/10666}.

\bibitem[Marquardt(1963)]{marquardtAlgorithmLeastSquares1963}
Donald~W. Marquardt.
\newblock An algorithm for least-squares estimation of nonlinear parameters.
\newblock \emph{Journal of the Society for Industrial and Applied Mathematics},
  11\penalty0 (2):\penalty0 431--441, 1963.
\newblock \doi{10.1137/0111030}.

\bibitem[Pal et~al.(2026)Pal, Holtorf, Larsson, Loman, Utkarsh, Schäfer, Qu,
  Edelman, and Rackauckas]{palNonlinearSolvejlHighPerformanceRobust2024}
Avik Pal, Flemming Holtorf, Axel Larsson, Torkel Loman, Utkarsh, Frank
  Schäfer, Qingyu Qu, Alan Edelman, and Christopher Rackauckas.
\newblock {NonlinearSolve}.jl: {High}-{Performance} and {Robust} {Solvers} for
  {Systems} of {Nonlinear} {Equations} in {Julia}.
\newblock \emph{ACM Transactions on Mathematical Software}, 52\penalty0
  (1):\penalty0 1--26, 2026.
\newblock \doi{10.1145/3779117}.

\bibitem[Ba{\c{s}}ar and Olsder(1999)]{basarDynamicNoncooperativeGame1999}
Tamer Ba{\c{s}}ar and Geert~Jan Olsder.
\newblock \emph{Dynamic Noncooperative Game Theory}, volume~23 of
  \emph{Classics in Applied Mathematics}.
\newblock Society for Industrial and Applied Mathematics, Philadelphia, PA, 2
  edition, 1999.
\newblock \doi{10.1137/1.9781611971132}.

\bibitem[Ab~Azar et~al.(2020)Ab~Azar, Shahmansoorian, and
  Davoudi]{abazarInverseOptimalControl2020}
Nematollah Ab~Azar, Aref Shahmansoorian, and Mohsen Davoudi.
\newblock From inverse optimal control to inverse reinforcement learning: A
  historical review.
\newblock \emph{Annual Reviews in Control}, 50:\penalty0 119--138, 2020.
\newblock \doi{10.1016/j.arcontrol.2020.06.001}.

\bibitem[Molloy et~al.(2022)Molloy, Inga~Charaja, Hohmann, and
  Perez]{molloyInverseOptimalControl2022}
Timothy~L. Molloy, Jairo Inga~Charaja, S{\"o}ren Hohmann, and Tristan Perez.
\newblock \emph{Inverse Optimal Control and Inverse Noncooperative Dynamic Game
  Theory: A Minimum-Principle Approach}.
\newblock Communications and Control Engineering. Springer, Cham, 2022.
\newblock \doi{10.1007/978-3-030-93317-3}.

\bibitem[Banga and Sager(2025)]{bangaGeneralizedInverseOptimal2025}
Julio~R. Banga and Sebastian Sager.
\newblock Generalized inverse optimal control and its application in biology.
\newblock \emph{Annual Reviews in Control}, 60:\penalty0 101029, 2025.
\newblock \doi{10.1016/j.arcontrol.2025.101029}.

\bibitem[Banga and Villaverde(2026)]{bangaHierarchicalFrameworkInverse2026}
Julio~R. Banga and Alejandro~F. Villaverde.
\newblock A hierarchical framework for inverse problems in biological
  cybernetics: Opportunities and limitations.
\newblock Preprint, DIGITAL.CSIC, 2026.

\bibitem[Martensen et~al.(2026)Martensen, Plate, Kirches, Potschka, Rackauckas,
  and Sager]{martensenetal_2026_corleonepaperexamples}
Carl~Julius Martensen, Christoph Plate, Christian Kirches, Andreas Potschka,
  Christopher~Vincent Rackauckas, and Sebastian Sager.
\newblock {Corleone.jl}: examples and benchmarks, 2026.
\newblock URL \url{https://github.com/mathopt/Corleone-Paper-Examples}.
\newblock Version v1.0.3.

\end{thebibliography}

\end{document}